\input epsf

\documentclass[12pt,a4paper]{article}
\usepackage{color}
\usepackage{colordvi}
\usepackage{amsmath,amssymb}
\usepackage{latexsym}
\usepackage{graphicx}
\usepackage{color}
\usepackage{multirow}
\usepackage{indentfirst}
\usepackage{subfigure}

\newtheorem{exam}{\hspace{6mm}Example}[section]

\begin{document}
\baselineskip=2pc

\centerline{\Large \bf A hybrid LDG-HWENO scheme for KdV-type equations%
\footnote{The research  was supported by NSFC grants 11328104 and 11571290.}}

\centerline{Dongmi Luo%
\footnote{School of Mathematical Sciences and Fujian Provincial
Key Laboratory of Mathematical Modeling and High-Performance
Scientific Computing, Xiamen University, Xiamen, Fujian 361005, China. E-mail: luodongmi@126.com.},
Weizhang Huang%
\footnote{Department of Mathematics, University of Kansas, Lawrence, Kansas 66045, U.S.A.
E-mail: whuang@ku.edu.},
and Jianxian Qiu%
\footnote{School of Mathematical Sciences and Fujian Provincial
Key Laboratory of Mathematical Modeling and High-Performance
Scientific Computing, Xiamen University, Xiamen, Fujian 361005, China. E-mail: jxqiu@xmu.edu.cn.}
}

\vspace{20pt}

\begin{abstract}

A hybrid LDG-HWENO scheme is proposed for the numerical solution of KdV-type partial differential equations.
It evolves the cell averages of the physical solution and its moments (a feature of Hermite WENO) while
discretizes high order spatial derivatives using the local DG method. The new scheme has the advantages
of both LDG and HWENO methods, including the ability to deal with high order spatial derivatives and
the use of a small number of global unknown variables. The latter is independent of
the order of the scheme and the spatial order of the underlying differential equations.
One and two dimensional numerical examples are presented to show
that the scheme can attain the same formal high order accuracy as the LDG method.
\end{abstract}


{\bf AMS(MOS) subject classification:} 65M60, 35L65
\pagenumbering{arabic}

\newpage

\section{Introduction}
\label{sec1}
\setcounter{equation}{0}
\setcounter{figure}{0}
\setcounter{table}{0}

We consider the numerical solution of KdV-type equations in one and two spatial dimensions.
The local discontinuous Galerkin (LDG)
method for this type of partial differential equations (PDEs) has been studied in \cite{EH05,E04}.
The LDG method is an extension of the discontinuous Galerkin (DG) method aimed at solving PDEs containing
higher than first order spatial derivatives.
The DG  method was first introduced by Reed and Hill \cite{EH08} for solving linear hyperbolic problems
for neutron transfer. A major development of the DG method was carried out by Cockburn et al. in a series of papers \cite{EH15,EH16,EH17,EH18}.

 The basic idea of the LDG method is to rewrite a PDE with high order spatial derivatives into a system of first order PDEs and
 then discretize it by the DG method. It can achieve nonlinear stability without slope limiters when carefully designed.
 The sub-optimal error estimates in $L^2$ norm
were obtained for the smooth solution of linear equations and the cell entropy inequality was
proven in \cite{E04}. While the LDG method also has many other advantages, it has the disadvantage of employing a large number of unknown variables and the number increases rapidly as the order of the method increases
especially in multiple dimensions and for high order PDEs. On the other hand, the essentially non-oscillatory (ENO) schemes of Shu and Osher \cite{EH13,EH14}, the weighted ENO (WENO) schemes of Liu et al. \cite{EH20} and
Jiang and Shu \cite{EH07}, and Hermite WENO (HWENO) schemes of Qiu and Shu \cite{CD02,BCCD03} have the advantage of employing a small number of unknown variables. But they have the disadvantage that the stencil used in reconstruction is becoming wider with an increasing order
of accuracy.

To avoid the disadvantages of DG/LDG and WENO schemes, based on the reconstruction procedure for HWENO limiter by Qiu and Shu \cite{CD02, BCCD03}  a hybrid RKDG-HWENO scheme was proposed by Balsara  et al. \cite{EH19} for hyperbolic conservation laws. It has intrinsic robustness and smaller stencils of the RKDG scheme and takes the advantage of the WENO scheme with a small number of unknown variables. The scheme evolves lower moments while reconstructing the higher moments of the solution. Dumbser et al. \cite{PnPm2008} extended the scheme to a new family of in-cell recovery DG method, referred to as $P_{N}P_{M}$ methods, where $P_{N}$ indicates that a piecewise polynomial of degree $N$ is used to represent a DG solution, and $P_{M}$ stands for a reconstructed polynomial solution of degree $M$ ($M \ge N$) that is used to compute the numerical fluxes. This approach yields a general, unified framework that
contains two important special cases, classical high order finite volume (FV) schemes $(N=0)$ and
the conventional discontinuous Galerkin (DG) method $(N=M)$, and has been applied successfully
to simulate the Euler equations of compressible gas dynamics and the equations of ideal
magnetohydrodynamics (MHD). $P_{N}P_{M}$ methods  were extended to the numerical solution
of the compressible Navier-Stokes equations \cite{PnPm2010, LLNMD}.

Motivated by  the hybrid RKDG-HWENO and $P_{N}P_{M}$ methods,  we propose a hybrid LDG-HWENO scheme
for the numerical solution of KdV-type PDEs. To our best knowledge, those methods have not been studied
for higher order PDEs like KdV-type equations containing third order spatial derivatives.
The new scheme employs LDG to approximate higher than first order spatial derivatives contained in the PDEs.
 Also different from HWENO methods in \cite{CD02,BCCD03}, where the cell averages of both the solution $u$ and its first derivative
$u_x$ are evolved in time and used in the reconstruction, we use the cell averages of $u$ (the solution) and
$u\frac{x-x_j}{\Delta x}$ (the first moment(s)) which are used in the reconstruction procedure for HWENO limiter \cite{CD02,BCCD03}.
The new scheme employs only a small number of (global) unknown variables and uses a small stencil
in the reconstruction. Numerical examples show that the new scheme can attain the same formal high order
accuracy as the LDG method.

An outline of the paper is given as follows. The hybrid LDG-HWENO scheme is described in Sections \ref{sec2}
and \ref{sec3} for one and two dimensional KdV-type equations, respectively.
In Section \ref{sec4}, a selection of one and two dimensional numerical examples are presented to demonstrate
the accuracy and the capability of the scheme. Conclusions are drawn in Section \ref{sec5}.

\section{The hybrid LDG-HWENO scheme in one dimension}
\label{sec2}
\setcounter{equation}{0}
\setcounter{figure}{0}
\setcounter{table}{0}

We consider the numerical solution of one dimensional KdV-type problems in the general form,
 \begin{equation}
 \label{u1}
 u_{t}+f(u)_{x}+(r'(u)g(r(u)_x)_x)_x=0 ,\quad a < x < b, \; t > 0
 \end{equation}
 with the initial condition
 \[
 u(x,0)=u_0(x),\quad a\leq x \leq b
 \]
 where $f(u), r(u)$, and $g(r)$ are given smooth functions. In this work, we consider a uniform mesh of cell size
 $\Delta x$. Denote the cells by $I_j=(x_{j-\frac{1}{2}},x_{j+\frac{1}{2}})$, where
 $x_{j+\frac{1}{2}}=\frac{1}{2}(x_j+x_{j+1})$. As for the HWENO scheme \cite{CD02,BCCD03}, we want
to find the approximations for the cell averages of $u$ and $u\frac{x-x_j}{\Delta x}$ of the solution of (\ref{u1}), i.e.,
 \[
 \Bar u_j \approx \frac{1}{\Delta x}\int _{I_j}udx,\quad
 \Bar v_j\approx \frac{1}{\Delta x}\int _{I_j}u\frac{x-x_j}{\Delta x}dx .
\]
Similar to the procedure of HWENO limiter for DG in  \cite{CD02,BCCD03}, we use here the first moment $\Bar v_j$ instead of the cell average of $u_x$ which was used for HWENO scheme  in  \cite{CD02,BCCD03}.
This is more consistent with the basis functions of LDG (see (\ref{basis-1}) below). We employ LDG (e.g., see \cite{EH05,E04}) for the discretization of high order spatial
 derivatives in (\ref{u1}). We first introduce the new variables
 \[
 q=r(u)_x,\qquad p=g(q)_x ,
 \]
 and rewrite (\ref{u1}) as a system of first order differential equations,
 \begin{align}
 \label{u2}
 u_{t}+(f(u)+r'(u)p)_x & =0, \\
 \label{u3}
 p & =g(q)_x, \\
 \label{u4}
 q &=r(u)_x .
 \end {align}
Next, we consider the discretization of (\ref{u3}) and (\ref{u4}).  Let $p$ and $q$ be  approximated by
\[
p_h =\sum\limits^k_{l=0}p^{(l)}_j\phi^{(j)}_l(x),\quad
q_h =\sum\limits^k_{l=0}q^{(l)}_j\phi^{(j)}_l(x), \quad \forall x \in I_j
\]
where $\phi^{(j)}_l(x)$'s are the orthogonal basis functions and $k$ is the highest degree of polynomials used in the basis functions. The first five basis functions are
\begin{align}
&\phi^{(j)}_0(x)=1,
&\phi^{(j)}_1(x)=\frac{x-x_j}{\Delta x},
\label{basis-1}\\
&\phi^{(j)}_2(x)=(\frac{x-x_j}{\Delta x})^2-\frac{1}{12},
& \phi^{(j)}_3(x)=(\frac{x-x_j}{\Delta x})^3-\frac{3}{20}\frac{x-x_j}{\Delta x},\notag \\
&\phi^{(j)}_4(x)=(\frac{x-x_j}{\Delta x})^4-\frac{3}{14}(\frac{x-x_j}{\Delta x})^2+\frac{3}{560} .&
\notag
\end{align}
Multiplying (\ref{u3}) and (\ref{u4}) by test functions $w$ and $z$, respectively,
integrating over the interval $I_j$, integrating by parts, and replacing function values at
cell boundaries by their numerical fluxes (those quantities with ``hats'') which we will define in the below,  we obtain
  \begin{align}
  \label{v1}
& \int_{I_{j}} p_h w dx+\int_{I_{j}} g(q_h)w_xdx-\Hat g_{j+\frac{1}{2}}w^-
_{j+\frac{1}{2}}+\Hat g_{j-\frac{1}{2}}w^+_{j-\frac{1}{2}}=0 ,
\\
\label{v2}
& \int_{I_{j}} q_h zdx+\int_{I_{j}} r(u) z_xdx-\Hat r_{j+\frac{1}{2}}z^-_{j+\frac{1}{2}}+\Hat r_{j-\frac{1}{2}}z^+_{j-\frac{1}{2}}=0,
\end {align}
where
\[
w_{j\pm \frac 1 2}^{\pm} = \lim_{x\to x_{j\pm \frac 1 2}\pm 0} w(x),
\quad z_{j\pm \frac 1 2}^{\pm} = \lim_{x \to x_{j\pm \frac 1 2}\pm 0} z(x).
\]
The test functions $w$ and $z$ in (\ref{v1}) and (\ref{v2}) are taken as the basis functions $\phi^{(j)}_l(x),\quad l=0,1,\cdots,k$, successively. The integrals in (\ref{v1}) and (\ref{v2}) are computed numerically using a Gaussian quadrature rule.  For the ($k+1$)th order accuracy, the integration must achieve $(2k+2)$th order as the DG method. Thus, we use the $(k+1)-$point Gaussian quadrature, namely,
$$\int_{I_j}g(q_h)w_xdx\approx\Delta x\sum_Gg(q_h(x_G))w_x(x_G)w_G,$$
$$\int_{I_j}r(u)z_xdx\approx\Delta x\sum_Gr(u(x_G))z_x(x_G)w_G,$$
where $w_G$'s are the weights and $x_G$'s are the Gauss points. The function $q_h$ used in the first integral is obtained by solving (\ref{v2}) while the value $u(x_G)$ used in the second integral is reconstructed from the cell averages \{$\Bar u_j,\Bar v_j$\}. For example, the sixth order reconstruction procedure includes following five steps.

Step 1.  Define the small stencils
 \[
 S_0=\{I_{j-1},I_j, I^\prime_{j-1}, I^\prime_j\},\quad
 S_1=\{I_{j},I_{j+1}, I^\prime_{j}, I^\prime_{j+1}\},\quad
 S_2=\{I_{j-1},I_j, I_{j+1}, I^\prime_j\}
 \]
 and a bigger stencil $S=\{S_0,S_1,S_2\}$, where $I_j$ and $I_j'$ stand for intervals chosen for
 $\Bar u_j$ and $\Bar v_j$, respectively. These stencils are used to  construct polynomials
 $p_0(x)$, $p_1(x)$, and $p_2(x)$ of degree three and polynomial $Q(x)$ of degree five,
 \begin{align*}
 & \frac{1}{\Delta x}\int_{I_{j+l}} p_0(x)dx=\Bar u_{j+l},
 \quad \frac{1}{\Delta x}\int_{I'_{j+l}=I_{j+l}} p_0(x)\frac{x-x_j}{\Delta x}dx=\Bar v_{j+l},\quad l=-1,0
 \\
 & \frac{1}{\Delta x}\int_{I_{j+l}} p_1(x)dx=\Bar u_{j+l},
 \quad \frac{1}{\Delta x}\int_{I'_{j+l}=I_{j+l}} p_1(x)\frac{x-x_j}{\Delta x}dx=\Bar v_{j+l},\quad l=0,1
 \\
 &\frac{1}{\Delta x}\int_{I_{j+l}} p_2(x)dx=\Bar u_{j+l},
 \quad \frac{1}{\Delta x}\int_{I'_{j}=I_{j}} p_2(x)\frac{x-x_j}{\Delta x}dx=\Bar v_{j},\quad l=-1,0,1
 \\
& \frac{1}{\Delta x}\int_{I_{j+l}} Q(x)dx=\Bar u_{j+l},
 \quad \frac{1}{\Delta x}\int_{I'_{j+l}=I_{j+l}} Q(x)\frac{x-x_j}{\Delta x}dx=\Bar v_{j+l},\quad l=-1,0,1.
 \end{align*}

Step 2. Compute the smoothness indicators, denoted by $\beta_l$ for each stencil $S_l$.
They measure how smooth $p_0(x),\; p_1(x)$, and $p_2(x)$ are in the target cell $I_j$.
Following \cite{EH07}, we define
\[
\beta_l=\sum\limits^3_{m=1}\int_{I_j}\Delta x^{2m-1}(\frac{d ^mp_l(x)}{d x^m})^2dx,\quad l=0,1,2.
\]

Step 3. For a given point $\Hat x \in I_j$, we find the linear weights, denoted by $\gamma_0,\gamma_1$ and $\gamma_2$, such that
\[
Q(\Hat x)=\sum\limits^2_{i=0}\gamma_i p_i(\Hat x).
\]
For example, for $\Hat x=x_{j+\frac{1}{2}}$ we have
\[
\gamma_0=\frac{25}{189},\quad \gamma_1=\frac{14}{27},\quad \gamma_2=\frac{22}{63}.
\]

Step 4. The nonlinear weights based on the smoothness indicators are computed as
\[
\omega_m=\frac{\Bar \omega_m}{\sum\limits_{l=0}^2\Bar \omega_l},\quad
\Bar \omega_l=\frac{\gamma_l}{(\lambda+\beta_l)^2},\quad m=0,1,2,
\]
where $\lambda$ is a small positive number to avoid the denominator to become zero.

Step 5.  Finally, the HWENO reconstruction is given by
\begin{equation}
u(\Hat x) = \sum\limits^2_{i=0}\omega_i p_i(\Hat x).
\label{u-0}
\end{equation}

Step 3 to Step 5 are repeated for each point $\Hat x \in \{\{x_G\},x_{j-\frac{1}{2}},x_{j+\frac{1}{2}}\}$, where $\{x_G\}$ denotes the set of the Gauss points.

 It is remarked that the reconstructed values of $u$ at both ends of $I_j$ are denoted by $u^+_{j-\frac{1}{2}}$ and $u^-_{j+\frac{1}{2}}$, respectively. Moreover, the linear weights for some points can become negative. For example, for $x_G=x_{j-\frac{0.5384693101056831}{2}}$,  the linear weights are $\gamma_0=-1.19876833424689$, $\gamma_1=-0.189130224626382$, and $\gamma_2=2.38789855887328$. The negative weights may lead to instability of the reconstruction. The technique developed by Shi et al.  \cite{EH06} can be used to treat reconstruction procedure with negative weights; the interested reader is referred to \cite{EH06} for the detail.

Multiplying (\ref{u2}) with $1$ and $\frac{x-x_j}{\Delta x}$ and integrating them by parts over $I_j$, we have
 \begin{align}
 \label{u5}
 \frac{d\Bar u_j}{dt}=-\frac{1}{\Delta x}(f(u)+r'(u)p)|_{I_j} ,
 \end {align}
 \begin{align}
 \label{u6}
 \frac{d\Bar v_j}{dt}=-\frac{1}{\Delta x}(f(u)+r'(u)p)\frac{x-x_j}{\Delta x}|_{I_j}
 +\frac{1}{\Delta x^2}\int_{I_j}(f(u)+r'(u)p) dx .
 \end {align}
The integral term in  (\ref{u6}) is approximated by the $(k+1)$-point Gaussian quadrature rule as for (\ref{v1}) and (\ref{v2}). We obtain the numerical scheme as
\begin{align}
\label{v3}
& \frac{d\Bar u_j}{dt}=-\frac{1}{\Delta x}[(\Hat f_{j+\frac{1}{2}}+\Hat r'_{j+\frac{1}{2}}\Hat p_{j+\frac{1}{2}})
-(\Hat f_{j-\frac{1}{2}}+\Hat r'_{j-\frac{1}{2}}\Hat p_{j-\frac{1}{2}})] ,
\\
\label{v4}
& \frac{d\Bar v_j}{dt}=-\frac{1}{2\Delta x}[(\Hat f_{j+\frac{1}{2}}+\Hat r'_{j+\frac{1}{2}}\Hat p_{j+\frac{1}{2}})
+(\Hat f_{j-\frac{1}{2}}+\Hat r'_{j-\frac{1}{2}}\Hat p_{j-\frac{1}{2}})]+\frac{1}{\Delta x^2}F_j ,
 \end {align}
 where \[F_j=\Delta x \sum_G(f(u(x_G))+r'(u(x_G))p_h(x_G))w_G\approx \int_{I_j}(f(u)+r'(u)p_h) dx.\]

A key component of the above described hybrid LDG-HWENO scheme is
to define the numerical fluxes in (\ref{v1}), (\ref{v2}), (\ref{v3}), and (\ref{v4}) to ensure the accuracy and stability of the scheme. We define
\begin{align*}
&\Hat f_{j+\frac{1}{2}}=\Hat f(u^-_{j+\frac{1}{2}},u^+_{j+\frac{1}{2}}),
&\Hat g_{j+\frac{1}{2}}=\Hat g(q^-_{j+\frac{1}{2}},q^+_{j+\frac{1}{2}}),\\
&\Hat p_{j+\frac{1}{2}}=p^+_{j+\frac{1}{2}},
&\Hat r_{j+\frac{1}{2}}=r(u^-_{j+\frac{1}{2}}),\\
&\Hat r'_{j+\frac{1}{2}}=\frac{r(u^+_{j+\frac{1}{2}})-r(u^-_{j+\frac{1}{2}})}{u^+_{j+\frac{1}{2}}-u^-_{j+\frac{1}{2}}} ,
\end{align*}
where
\[
p_{j+\frac12}^{\pm} = p_h(x_{j+\frac12}\pm 0),\quad q_{j+\frac12}^{\pm} = q_h(x_{j+\frac12}\pm 0),
\]
and $u_{j+\frac 1 2}^{\pm}$ are the reconstructed values of $u$  by the HWENO method as described in Step~1-Step~5. It is noted that the choice of the numerical fluxes is not unique. The key part is that $\Hat p$ and $\Hat r$ must be taken from the opposite sides. The numerical flux $\Hat f(a,b)$ for $f(u)$ should satisfy the following conditions:
\begin{enumerate}
\item[(i)] $\Hat f(a,b)$ is a Lipschitz continuous function in both arguments $a$ and $b$;
\item[(ii)] $\Hat f$ is consistent with $f(u)$, namely, $\hat f(u,u)=f(u)$;
\item[(iii)] $\Hat f(a,b)$ is a monotone flux, i.e., it is non-decreasing in $a$ and non-increasing in  $b$.
\end{enumerate}
In this paper, we use the local Lax-Friedrichs(LLF) flux,
\[
\Hat f(a,b)=\frac{1}{2}[f(a)+f(b)-\alpha (b-a)], \quad \alpha=\max\limits_{u\in [\min(a,b),\max(a,b)]}|f'(u)| .
\]

The numerical flux $\Hat g(c,d)$ for $g(q)$ should also satisfy conditions (i) and
(ii) but the condition (iii) should be replaced by
\begin{enumerate}
\item[(iii)'] $-\Hat g(q^-,q^+)$ is a monotone flux for $-g(q)$, namely, $\Hat g(c,d)$ is non-increasing in $c$
and non-decreasing in  $d$.
\end{enumerate}
Corresponding to the LLF flux, we have
\[
\Hat g(c,d)=\frac{1}{2}[g(c)+g(d)-\alpha (c-d)] ,\quad \alpha=\max_{q\in [\min(c,d),\max(c,d)]}|g'(q)| .
\]

We note that (\ref{v2}) can be solved independently on each cell for $q_h$
in terms of the cell averages $\{\Bar u_j,\Bar v_j\}$. Similarly, (\ref{v1}) can be solved on each cell
for $p_h$ in terms of $q_h$. As a consequence, the global unknown variables involve only
$\Bar u_j $ and $\Bar v_j$, which are governed by a system of ordinary differential equations.
In our computation, the system is integrated in time using an explicit TVD Runge-Kutta method,  such as the 3rd TVD Runge-Kutta method for solving
$$u_t=L(u,t),$$
where $L(u,t)$ is a spatial discretization operator,
\begin{align}
&u^*=u^n+\Delta tL(u^n,t^n) \notag\\
&u^{**}=\frac{3}{4}u^n+\frac{1}{4}(u^*+\Delta tL(u^*,t^n+\Delta t)) \notag \\
&u^{n+1}=\frac{1}{3}u^n+\frac{2}{3}(u^{**}+\Delta tL(u^{**},t^n+\frac{1}{2}\Delta t)). \notag
\end{align}

From the construction of the scheme, it is not difficult to show formally that the method is of $(k+1)^{\text{th}}$ order
in space ($k \le 4$) and third order in time, i.e.,
\[
e_h^n = \mathcal{O}(\Delta t^3) + \mathcal{O}(\Delta x^{k+1}).
\]
Since an explicit scheme is used, the time step is subject to the CFL condition and should be
taken as $\Delta t=\mathcal{O}(\Delta x^3)$. For this choice, the error is dominated by the spatial
discretization error.

\section{The hybrid LDG-HWENO method in two dimensions}
\label{sec3}
\setcounter{equation}{0}
\setcounter{figure}{0}
\setcounter{table}{0}

In this section we describe the method for two dimensional problems.
We consider the general form
 \begin{align}
 &u_{t}+f_1(u)_{x}+f_2(u)_y+[r'_1(u)(g_{11}(r_1(u)_x)_x+g_{12}(r_1(u)_x)_y)]_x \notag\\
 & +[r'_2(u)(g_{21}(r_2(u)_y)_x+g_{22}(r_2(u)_y)_y)]_y=0,
 \quad (x, y) \in (a,b)\times (c,d)
 \label{u2d-1}
 \end {align}
subject to a periodic boundary condition and the initial condition
 \[
 u(x,y,0)=u_0(x,y),\quad (x, y) \in (a,b)\times (c,d)
\]
where $f_m(u)$, $r_m(u)$, and $g_{mn}(r)$ are given smooth nonlinear functions.

We use a rectangle mesh of cell size  $\Delta x$ and $\Delta y$ in $x$ and $y$ directions,
respectively. We denote the cells by
\[
I_{ij}=(x_{i-\frac{1}{2}},\, x_{i+\frac{1}{2}})
\times (y_{j-\frac{1}{2}},\, y_{j+\frac{1}{2}}),
\]
where
\[
x_{i+\frac{1}{2}} = \frac{1}{2}(x_{i}+x_{i+1}),\quad
y_{j+\frac{1}{2}} = \frac{1}{2}(y_{j}+y_{j+1}) .
\]
We approximate the cell averages of $u$, $u\frac{x-x_{i}}{\Delta x}$,  $u \frac{y-y_{j}}{\Delta y}$ and $u \frac{x-x_{i}}{\Delta x} \frac{y-y_{j}}{\Delta y}$ by
\begin{align*}
& \Bar u_{ij}\approx \frac{1}{\Delta x\Delta y}\int _{I_{ij}}udxdy,\\
& \Bar v_{ij}\approx \frac{1}{\Delta x\Delta y}\int _{I_{ij}}u\frac{x-x_{i}}{\Delta x}dxdy,\\
&\Bar w_{ij}\approx \frac{1}{\Delta x\Delta y}\int _{I_{ij}}u\frac{y-y_{j}}{\Delta y}dxdy ,\\
&\Bar Z_{ij}\approx \frac{1}{\Delta x\Delta y}\int _{I_{ij}}u\frac{x-x_{i}}{\Delta x}\frac{y-y_{j}}{\Delta y}dxdy .
\end{align*}

For the discretization of high order spatial derivatives using LDG, we introduce the new variables
 \begin{equation}
 q_1=r_1(u)_x,\quad q_2=r_2(u)_y,\quad p_1=g_{11}(q_1)_x+g_{12}(q_1)_y,\quad p_2=g_{21}(q_2)_x+g_{22}(q_2)_y,
 \label{uuu}
 \end{equation}
  and rewrite (\ref{u2d-1}) as a first order system
  \begin {align}
  \label{uu1}
  u_t+(f_1(u)+r'_1(u)p_1)_{x}+(f_2(u)+r'_2(u)p_2)_{y}=0 ,
  \end {align}
coupled with (\ref{uuu}).
We approximate $p_m,q_m(m=1,2)$ by polynomials,
$$p_m=\sum\limits_{l=0}^kp_{m,l}^{(ij)}\phi^{(ij)}_l(x,y),\quad q_m=\sum\limits_{l=0}^kq_{m,l}^{(ij)}\phi^{(ij)}_l(x,y) ,\quad m=1,2\quad \forall (x,y)\in I_{ij}$$
where $\phi^{(ij)}_l(x,y)$'s are the basis functions,
$$\phi^{(ij)}_0=1,\quad \phi^{(ij)}_1=\frac{x-x_{i}}{\Delta x},\quad\phi^{(ij)}_2=\frac{y-y_{j}}{\Delta y},\quad\phi^{(ij)}_3=(\frac{x-x_{i}}{\Delta x})^2-\frac{1}{12}$$
$$\phi^{(ij)}_4=(\frac{x-x_{i}}{\Delta x})(\frac{y-y_{j}}{\Delta y}),\quad\phi^{(ij)}_5=(\frac{y-y_{j}}{\Delta y})^2-\frac{1}{12},\quad\cdots$$

To discretize (\ref{uuu}), we multiply it by test functions $w_m$ and $z_m$, respectively,
and integrate the resulting equations by parts over the cell $I_{ij}$. Using numerical fluxes for quantities on the boundary of the cell, we have
  \begin {align}
  \label{uu2}
 & \int_{I_{ij}}p_1w_1dxdy+\int_{I_{ij}}g_{11}(q_1)(w_1)_{x}dxdy +\int_{I_{ij}}g_{12}(q_1)(w_1)_{y}dxdy
  \notag \\
&  - \int_{y_{j-\frac{1}{2}}}^{y_{j+\frac{1}{2}}}\left [\Hat g_{11}(q_1(x_{i+\frac{1}{2}},y))w_1(x^-_{i+\frac{1}{2}},y)
-  \Hat g_{11}(q_1(x_{i-\frac{1}{2}},y))w_1(x^+_{i-\frac{1}{2}},y)\right ] dy\notag\\
 & - \int_{x_{i-\frac{1}{2}}}^{x_{i+\frac{1}{2}}}\left [ \Hat g_{12}(q_1(x,y_{j+\frac{1}{2}}))w_1(x,y^-_{j+\frac{1}{2}})
  -\Hat g_{12}(q_1(x,y_{j-\frac{1}{2}}))w_1(x,y^+_{j-\frac{1}{2}})\right ] dx=0,
  \end {align}

  \begin {align}
  \label{uuu2}
  & \int_{I_{ij}}p_2w_2dxdy+\int_{I_{ij}}g_{21}(q_2)(w_2)_{x}dxdy +\int_{I_{ij}}g_{22}(q_2)(w_2)_{y}dxdy
  \notag \\
  & - \int_{y_{j-\frac{1}{2}}}^{y_{j+\frac{1}{2}}}\left [\Hat g_{21}(q_2(x_{i+\frac{1}{2}},y))w_2(x^-_{i+\frac{1}{2}},y)
-  \Hat g_{21}(q_2(x_{i-\frac{1}{2}},y))w_2(x^+_{i-\frac{1}{2}},y)\right ] dy\notag\\
& -  \int_{x_{i-\frac{1}{2}}}^{x_{i+\frac{1}{2}}}\left [ \Hat g_{22}(q_2(x,y_{j+\frac{1}{2}}))w_2(x,y^-_{j+\frac{1}{2}})-
  \Hat g_{22}(q_2(x,y_{j-\frac{1}{2}}))w_2(x,y^+_{j-\frac{1}{2}})\right ] dx=0,
  \end {align}

  \begin{align}
  \label{uuu3}
&  \int_{I_{ij}}q_1z_1dxdy+\int_{I_{ij}}r_{1}(u)(z_1)_{x}dxdy \notag \\
&  -\int_{y_{j-\frac{1}{2}}}^{y_{j+\frac{1}{2}}}\left [ \Hat r_1(u(x_{i+\frac{1}{2}},y))z_1(x^-_{i+\frac{1}{2}},y)
- \Hat r_1(u(x_{i-\frac{1}{2}},y))z_1(x^+_{i-\frac{1}{2}},y)\right ]dy=0,
\end{align}

  \begin {align}
  \label{uuu4}
&  \int_{I_{ij}}q_2z_2dxdy+\int_{I_{ij}}r_{2}(u)(z_2)_{y}dxdy\notag \\
&  - \int_{x_{i-\frac{1}{2}}}^{x_{i+\frac{1}{2}}}\left [ \Hat r_2(u(x,y_{j+\frac{1}{2}}))z_2(x,y^-_{j+\frac{1}{2}})
-\Hat r_2(u(x,y_{j-\frac{1}{2}}))z_2(x,y^+_{j-\frac{1}{2}})\right ] dx=0 .
\end {align}
This leads to equations including double integrals and line integrals on the cell $I_{ij}$. As for the one dimensional case, the volume integrals can be approximated by the tensor product Gaussian quadrature rule with $(k+1)$ points in each direction, and the line integrals can be computed by the $(k+1)$-point Gaussian quadrature rule. 

The final discrete equations are obtained by replacing those the volume and line integrals
by the Gaussian quadrature formulas. To save space, we omit these equations here. Instead, we point out that
$p_1,\, p_2,\, q_1,\, q_2$ are obtained by solving (\ref{uu2})-(\ref{uuu4})  and the values of $u$ at $(\Hat x,\Hat y)$ with $\Hat x\in\{x_{i-\frac{1}{2}}, x_{i+\frac{1}{2}}, \{x_{G_1}\}\}$ and $\Hat y\in\{y_{j-\frac{1}{2}}, y_{j+\frac{1}{2}}, \{y_{G_2}\}\}$ are reconstructed from the cell averages as in the one dimensional case. For reconstruction on Cartesian meshes, one can adopt either a direct two dimensional procedure or a dimension-by-dimension strategy \cite{EH06}.  In this paper, we use the dimension-by-dimension strategy. First, we perform two $y$-direction reconstructions, i.e.,
\[
\{\Bar u_{mn},\Bar w_{mn}\}\longrightarrow \Bar u_{i+l,j}(\Hat y)\approx \frac{1}{\Delta x}\int_{I_{i+l,j}}u(x,\Hat y)dx,\quad l=-1,0,1,\quad \Hat y\in\{y_{j-\frac{1}{2}}, y_{j+\frac{1}{2}}, \{y_{G_2}\}\}
\]
\[
\{\Bar v_{mn},\Bar Z_{mn}\}\longrightarrow \Bar u_{x,i+l,j}(\Hat y)\approx \frac{1}{\Delta x}\int_{I_{i+l,j}}u_x(x,\Hat y)dx,\quad l=-1,0,1,\quad \Hat y\in\{y_{j-\frac{1}{2}}, y_{j+\frac{1}{2}}, \{y_{G_2}\}\}.
\]
Then we use $\Bar u(\Hat y)$ and $\Bar u_x(\Hat y)$ to perform $x$-direction reconstruction to get an approximation to $u(\Hat x,\Hat y)$, i.e.,
\[
\{\Bar u_{mn}(\Hat y),\Bar u_{x,m,n}(\Hat y)\}\longrightarrow \widetilde{u}(\Hat x,\Hat y)\approx u(\Hat x,\Hat y),\quad \Hat x\in\{x_{i-\frac{1}{2}}, x_{i+\frac{1}{2}}, \{x_{G_1}\}\},\quad \Hat y\in\{y_{j-\frac{1}{2}}, y_{j+\frac{1}{2}}, \{y_{G_2}\}\} .
\]
The values $\widetilde{u}(\Hat x,\Hat y)$ are used in computing the volume and line integrals.


Multiplying (\ref{uu1}) with $1$, $\frac{x-x_{i}}{\Delta x}$,
$\frac{y-y_{j}}{\Delta y}$ and $\frac{x-x_{i}}{\Delta x}\frac{y-y_{j}}{\Delta y}$, integrating over $I_{ij}$  by parts, and using numerical fluxes on the cell boundary, we get
\begin {align}
\frac{d\Bar u_{ij}}{dt} = & -\frac{1}{\Delta x\Delta y}\int_{y_{j-\frac{1}{2}}}^{y_{j+\frac{1}{2}}}\Hat H_1(x_{i+\frac{1}{2}},y)-\Hat H_1(x_{i-\frac{1}{2}},y)dy\notag\\
&  -\frac{1}{\Delta x\Delta y}\int_{x_{i-\frac{1}{2}}}^{x_{i+\frac{1}{2}}}\Hat H_2(x,y_{j+\frac{1}{2}})-\Hat H_2(x,y_{j-\frac{1}{2}})dx,
\end {align}

\begin {align}
\frac{d\Bar v_{ij}}{dt}= & -\frac{1}{2\Delta x\Delta y}\int_{y_{j-\frac{1}{2}}}^{y_{j+\frac{1}{2}}}\Hat H_1(x_{i+\frac{1}{2}},y)+\Hat H_1(x_{i-\frac{1}{2}},y)dy \notag\\
& -\frac{1}{\Delta x\Delta y}\int_{x_{i-\frac{1}{2}}}^{x_{i+\frac{1}{2}}}(\Hat H_2(x,y_{j+\frac{1}{2}})-\Hat H_2(x,y_{j-\frac{1}{2}}))\frac{x-x_{i}}{\Delta x}dx \notag \\
& + \frac{1}{\Delta x^2\Delta y}\int_{I_{ij}}H_1dxdy ,
\end {align}

\begin {align}
\label{uu3}
\frac{d\Bar w_{ij}}{dt}=& -\frac{1}{\Delta x\Delta y}\int_{y_{j-\frac{1}{2}}}^{y_{j+\frac{1}{2}}}(\Hat H_1(x_{i+\frac{1}{2}},y)-\Hat H_1(x_{i-\frac{1}{2}},y))\frac{y-y_{j}}{\Delta y}dy \notag\\
& -\frac{1}{2\Delta x\Delta y}\int_{x_{i-\frac{1}{2}}}^{x_{i+\frac{1}{2}}}(\Hat H_2(x,y_{j+\frac{1}{2}})+\Hat H_2(x,y_{j-\frac{1}{2}}))dx \notag \\
& + \frac{1}{\Delta x\Delta y^2}\int_{I_{ij}}H_2dxdy,
\end {align}

 \begin {align}
 \label{zz}
\frac{d\Bar Z_{ij}}{dt}= & -\frac{1}{2\Delta x\Delta y}\int_{y_{j-\frac{1}{2}}}^{y_{j+\frac{1}{2}}}(\Hat H_1(x_{i+\frac{1}{2}},y)+\Hat H_1(x_{i-\frac{1}{2}},y))\frac{y-y_{j}}{\Delta y}dy \notag \\
& -\frac{1}{2\Delta x\Delta y}\int_{x_{i-\frac{1}{2}}}^{x_{i+\frac{1}{2}}}(\Hat H_2(x,y_{j+\frac{1}{2}})+\Hat H_2(x,y_{j-\frac{1}{2}}))\frac{x-x_{i}}{\Delta x}dx \notag \\
& + \frac{1}{\Delta x^2\Delta y}\int_{I_{ij}}H_1\frac{y-y_{j}}{\Delta y}dxdy+\frac{1}{\Delta x\Delta y^2}\int_{I_{ij}}H_2\frac{x-x_{i}}{\Delta x}dxdy,
\end {align}
where $H_1=f_1(u)+r'_1(u)p_1$ and $H_2=f_2(u)+r'_2(u)p_2$. 

The volume integrals and line integrals are approximated  by a Gaussian quadrature rule as for (\ref{uu2})-(\ref{uuu4}).
If we use $G_b$ to stand for a boundary point $(\Hat x,\Hat y)$, then the numerical fluxes $\Hat H_i(G_b)$, $\Hat r_i(G_b)$, and $\Hat g_{ij}(G_b)$ $(i=1,2,j=1,2)$ can be defined similarly as in one dimension, viz.,
\begin{align*}
&\Hat f_i(u(G_b))=\Hat f_i(u^-(G_b),u^+(G_b)),& \Hat g_{ij}=\Hat g_{ij}(q_i^-(G_b),q_i^+(G_b)),
\\
&\Hat p_i=p_i^+(G_b),& \Hat r_i=r_i(u^-(G_b)),
\\
&\Hat r'_i(u(G_b))=\frac{r_i(u^+(G_b))-r_i(u^-(G_b))}{u^+(G_b)-u^-(G_b)},& \Hat H_i=\Hat f_i+\Hat r'_i\Hat p_i,
\end{align*}
where $u^{\pm}(G_b)$, $q_i^{\pm}(G_b)$, and $p_i^{\pm}(G_b)$ are the left (or ``in'') and right (or ``out'') limits of the solutions
$u$, $q_i$, and $p_i$ at the cell interface $G_b$, respectively.
As in one dimension, $\Hat f_i(u(G_b))$ is a monotone flux for $\Hat f_i(u,u)=f_i(u)$
and $-\Hat g_{ij}$ is a flux for $-\Hat g_{ij}(q_i,q_i)=-g_{ij}(q_i)$. Also, the choice of the fluxes is not unique. We must take $\Hat p$ and $\Hat r$ from the opposite sides.

Once again, the discrete equations can be obtained accordingly. To save space, they are omitted here.

The resultant ODE system is integrated in time with an explicit third order TVD Runge-Kutta method.

\section{Numerical examples}
\label{sec4}
\setcounter{equation}{0}
\setcounter{figure}{0}
\setcounter{table}{0}

In this section we present numerical results obtained with the hybrid LDG-HWENO scheme described
in the previous sections for four examples each in one and two dimensions.

\subsection{One dimensional examples}

\begin{exam}{\em
\label{exam4.1}

We compute the solution of the linear equation
\[
u_t+ u_{xxx}=0,\quad 0 < x < 2 \pi,\; t > 0
\]
subject to the initial condition $u(x,0)=\hbox {sin}(x)$ and the periodic boundary condition. The equation is in the form of (\ref{u1}) with $f(u)=0,\quad r(u)=u,\quad g(r)=r$. The exact solution is given by $u(x,t)=\hbox {sin}(x+ t)$. The $L^1$, $L^2$, and $L^{\infty}$ norm of the error and the convergence order are shown in Table \ref{ex1.2}. It can be seen that the scheme with $P^k$ elements in DG gives at least a $(k+1)^{\text{th}}$ order of accuracy.
}\end{exam}

\begin{table}
\caption{Example~\ref{exam4.1}: periodic boundary conditions, and $T=1.0$.}
\renewcommand{\multirowsetup}{\centering}
\begin{center}
\begin{tabular}{|c|c|c|c|c|c|c|}
\hline
k &n & 10 & 20 & 40 & 80 & 160\\
\hline
\multirow{6}{1cm}{2}
 & $L^1$       & 2.668e-3 & 1.830e-4 & 1.904e-5 & 2.285e-6 & 2.833e-7  \\
 &  Order      & \quad    & 3.87     & 3.26     & 3.06      &  3.01    \\
 &$L^2$       & 2.751e-3 & 1.959e-4 & 2.078e-5 & 2.516e-6 &  3.135e-7  \\
 &Order       & \quad    & 3.81     & 3.24     & 3.05      &   3.00    \\
 &$L_{\infty}$&3.593e-3  &2.571e-4  & 2.854e-5 & 3.522e-6  &  4.412e-7 \\
 &Order    & \quad       & 3.80     & 3.17     & 3.02      &   3.00    \\
 \hline
 \multirow{6}{1cm}{3}

 & $L^1$       & 2.772e-4 & 6.544e-6 & 1.805e7 & 9.452e-9 &  5.755e-10 \\
 &  Order      & \quad    & 5.40     & 5.18     & 4.26     & 4.04      \\
 &$L^2$       & 2.973e-4 & 6.886e-6 & 2.102e-7 & 1.063e-8 &  6.412e-10 \\
 &Order       & \quad    & 5.43     & 5.03     & 4.31       & 4.05      \\
 &$L_{\infty}$&4.321e-4  &9.798e-6  & 3.305e-7 & 1.605e-8  & 9.207e-10 \\
 &Order    & \quad       & 5.46     & 4.89     & 4.36      & 4.12      \\
 \hline
\multirow{6}{1cm}{4}
 &$L^1$       & 9.220e-4 & 2.399e-5 & 5.009e-7 & 1.037e-8  & 2.216e-10 \\

 &Order      & \quad    & 5.26     & 5.58     & 5.59      & 5.55 \\
 &$L^2$       & 1.139e-3 & 2.545e-5 & 5.251e-7 & 1.097e-8 & 2.387e-10 \\
 &Order       & \quad    & 5.48     & 5.60     & 5.58      & 5.53  \\
 &$L_{\infty}$&1.655e-3  &3.416e-5  & 6.987e-7 & 1.657e-8  & 3.762e-10 \\
 &Order    & \quad       & 5.60     & 5.61     & 5.40      & 5.46  \\
\hline
\end{tabular}
\end{center}
\label{ex1.2}
\end{table}

\vspace{20pt}

\begin{exam}{\em
\label{exam4.2}
In order to see the accuracy of the scheme for nonlinear problems, we compute the classical soliton
solution of the KdV equation
\[
u_t-3(u^2)_x+u_{xxx}=0,\quad -10\le x\le 12 .
\]
The initial condition is given by
\[
u(x,0)=-2\hbox {sech}^2(x),
\]
and the exact solution is
\[
u(x,t)=-2\hbox {sech}^2(x-4t) .
\]
For this example, $f(u)=-3u^2,\quad r(u)=u,\quad g(r)=r$. Table \ref{ex2} gives the error of the numerical solution at $t=0.5$. We can see that
the $(k+1)^{\text{th}}$ order of accuracy of the scheme
is achieved for this nonlinear problem.
}\end{exam}

\begin{table}
\caption{Example~\ref{exam4.2}: periodic boundary conditions and $T=0.5$.}

\renewcommand{\multirowsetup}{\centering}
\begin{center}
\begin{tabular}{|c|c|c|c|c|c|c|}
\hline
k &n           & 40 & 80 & 160 & 320 & 640\\
\hline
\multirow{6}{1cm}{2}

 & $L^1$       & 1.303e-2 & 1.126e-3 & 9.841e-5 & 1.103e-5  & 1.205e-6      \\
 &  Order      & \quad    & 3.53     & 3.52     & 3.28      &  3.07     \\
 &$L^2$       & 2.569e-2 & 1.937e-3 & 1.514e-4 & 1.548e-5   & 1.843e-6  \\
 &Order       & \quad    & 3.73     & 3.68     & 3.29       &  3.07     \\
 &$L_{\infty}$&7.996e-2  &8.245e-3  & 5.926e-4 & 5.992e-6  &  7.115e-6  \\
 &Order    & \quad       & 3.28     & 3.80     & 3.31      &   3.07     \\
\hline
\multirow{6}{1cm}{3}
 & $L^1$       & 9.391e-3 & 1.527e-4 & 5.075e-6 & 1.586e-7 & 7.290e-9 \\
 &  Order      & \quad    & 5.94     & 4.91     & 5.00      &  4.44\\
 &$L^2$       & 1.154e-2 & 2.443e-4 & 8.117e-6 & 2.567e-7 &1.169e-8  \\
 &Order       & \quad    & 5.56     & 4.91     & 4.98      & 4.46  \\
 &$L_{\infty}$&2.425e-2  &1.244e-3  & 3.427e-5 & 1.002e-6  & 4.141e-8 \\
 &Order    & \quad       & 4.28     & 5.18     & 5.10      &   4.60\\
\hline
\multirow{6}{1cm}{4}
 &$L^1$       & 2.634e-2 & 1.986e-3 & 1.402e-5 & 3.708e-7  & 1.024e-8 \\
 &Order      & \quad    & 3.73     & 7.15     & 5.24      & 5.18 \\
 &$L^2$       & 6.395e-2 & 3.346e-3 & 2.247e-5 & 5.747e-7 & 1.656e-8 \\
 &Order       & \quad    & 4.26     & 7.22     & 5.29      & 5.12  \\
 &$L_{\infty}$&2.906e-1  &1.197e-2  & 1.001e-4 & 2.587e-6  & 7.374e-8 \\
 &Order    & \quad       & 4.60     & 6.90     & 5.27      & 5.13  \\
\hline
\end{tabular}
\end{center}
\label{ex2}
\end{table}

\begin{figure}[hbtp]
  \begin{center}
  \mbox{\subfigure[t=1]
  {\includegraphics[width=4.5cm]{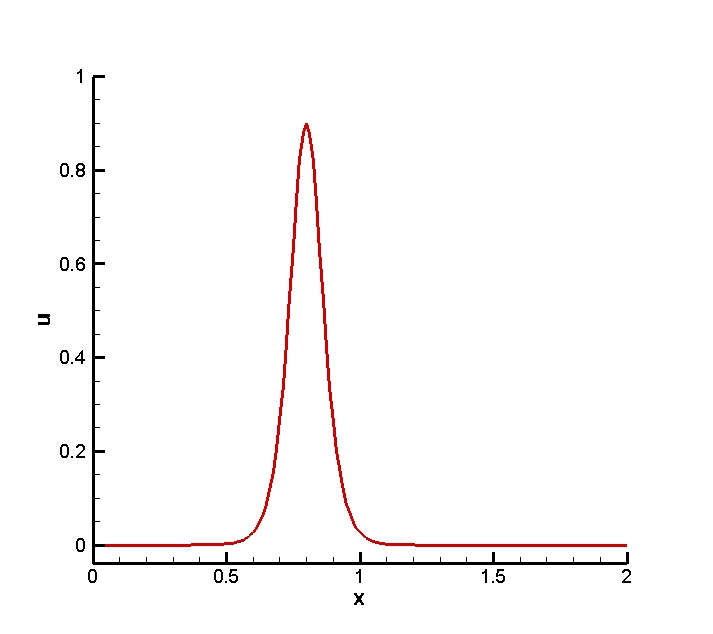}}\quad
    \subfigure[t=2]
    {\includegraphics[width=4.5cm]{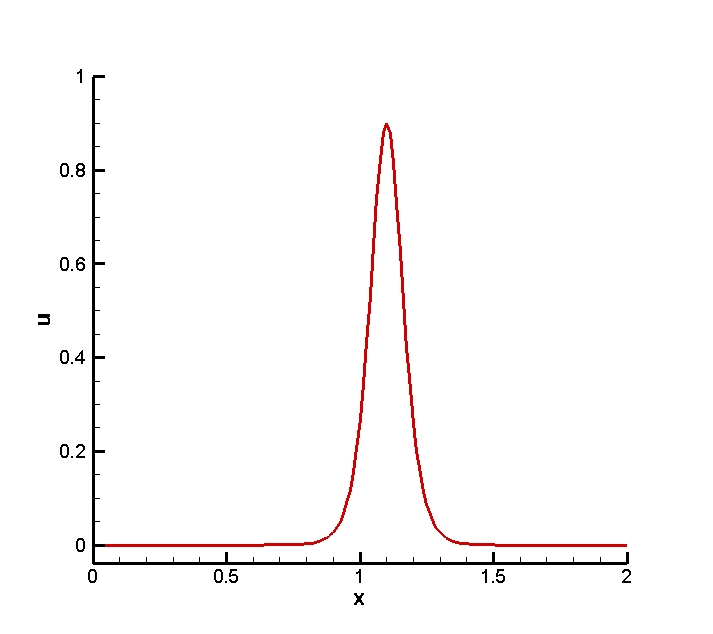}}
    }
  \mbox{\subfigure[t=1]
  {\includegraphics[width=4.5cm]{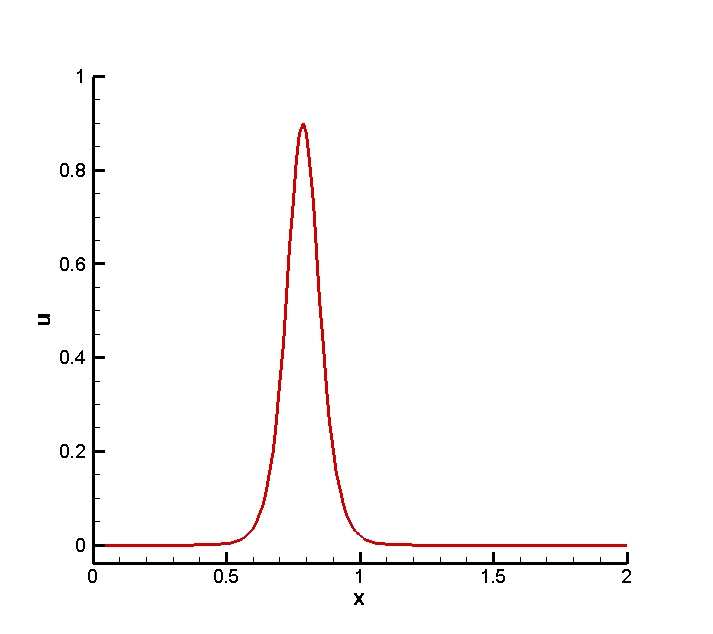}}\quad
    \subfigure[t=2]
    {\includegraphics[width=4.5cm]{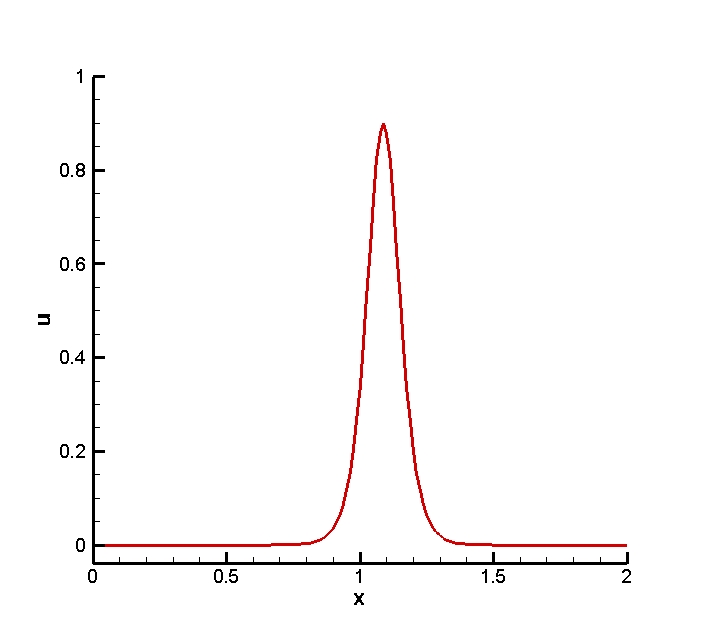}}
    }

    \caption{Example~\ref{exam4.3}: Single soliton case. Top: $P^2$ elements are used with 160 cells; Bottom: $P^4$ elements are used with 160 cells.}
    \label{fig:edge1}
    \end{center}
    \end{figure}

\begin{figure}[hbtp]
  \begin{center}
  \mbox{\subfigure[t=1]
  {\includegraphics[width=4.5cm]{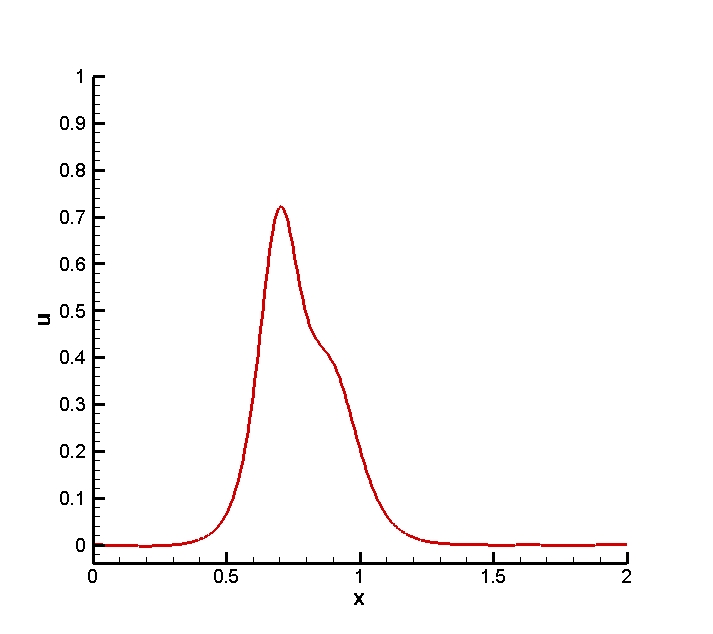}}\quad
    \subfigure[t=2]
    {\includegraphics[width=4.5cm]{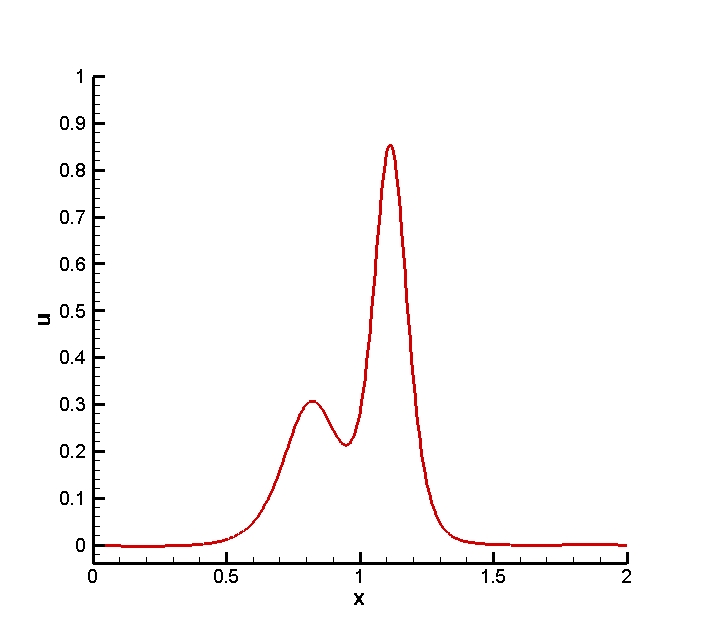}}
    }
    \mbox{\subfigure[t=1]
  {\includegraphics[width=4.5cm]{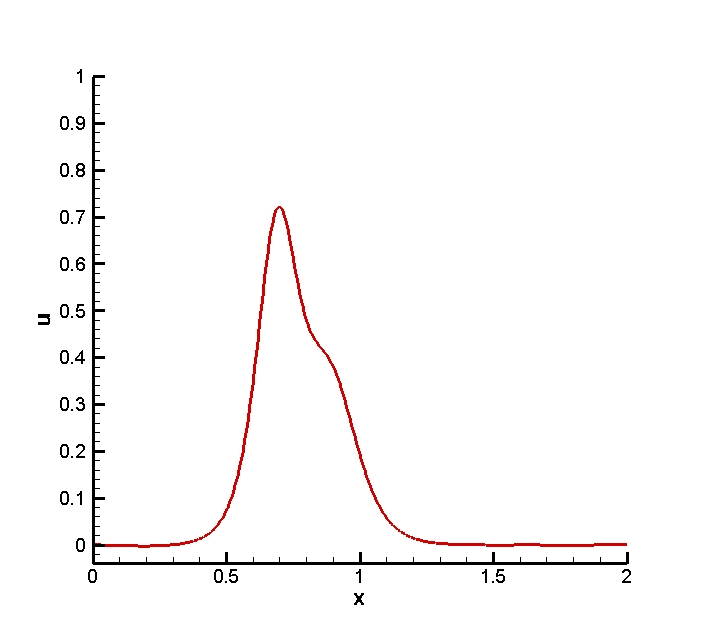}}\quad
    \subfigure[t=2]
    {\includegraphics[width=4.5cm]{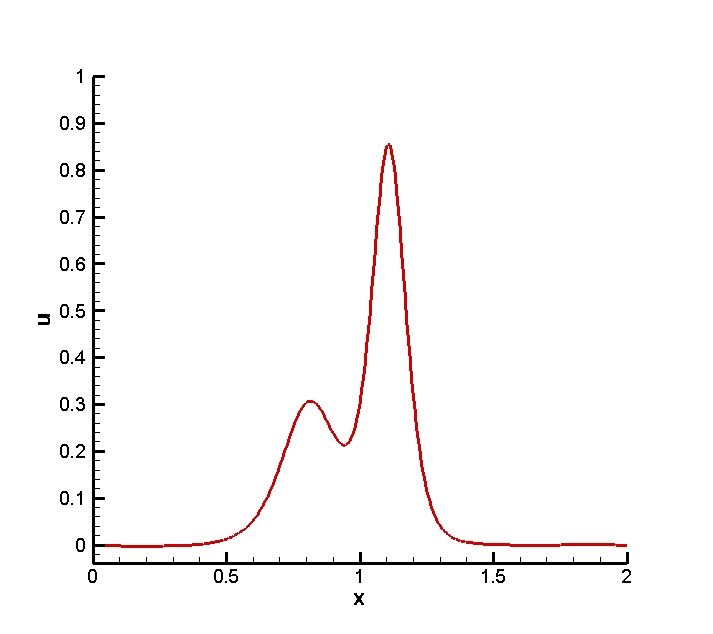}}
    }
    \caption{Example~\ref{exam4.3}: Double soliton collision case. Top: $P^2$ elements are used with 320 cells; Bottom: $P^4$ elements are used with 320 cells.}
    \label{fig:edge2}
    \end{center}
    \end{figure}

\begin{figure}[hbtp]
  \begin{center}
  \mbox{\subfigure[t=1]
  {\includegraphics[width=4.5cm]{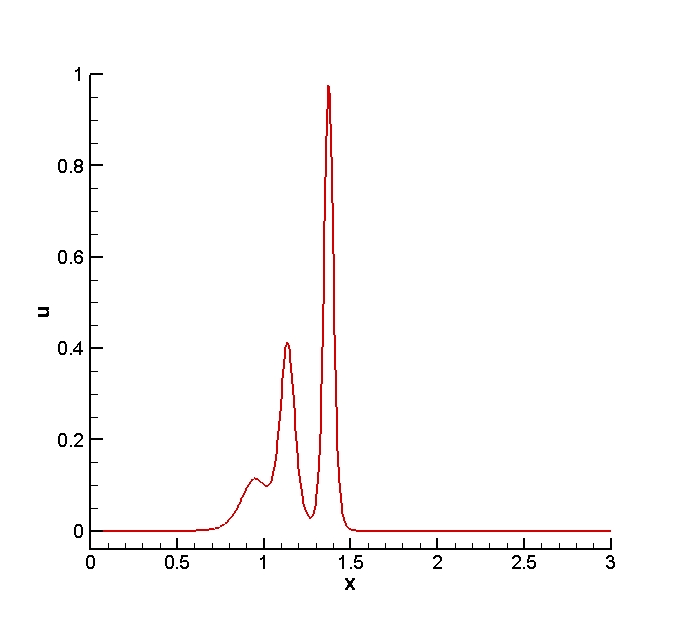}}\quad
    \subfigure[t=2]
    {\includegraphics[width=4.5cm]{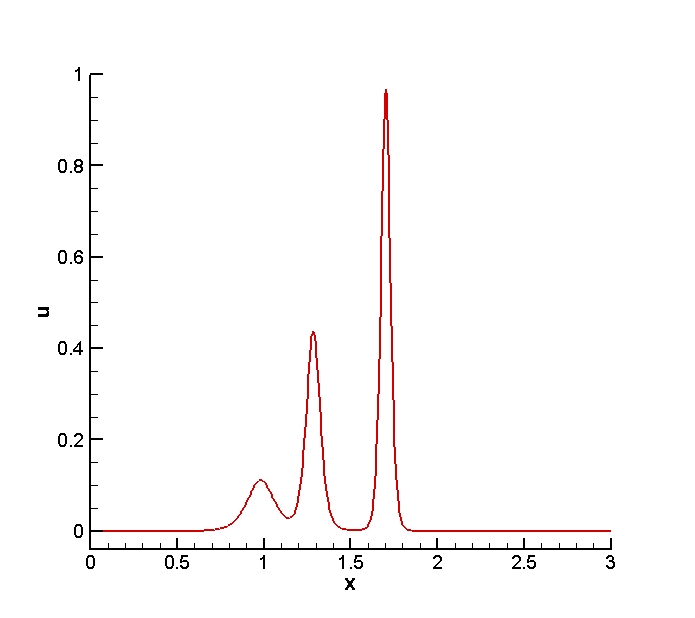}}
    }
  \mbox{\subfigure[t=1]
  {\includegraphics[width=4.5cm]{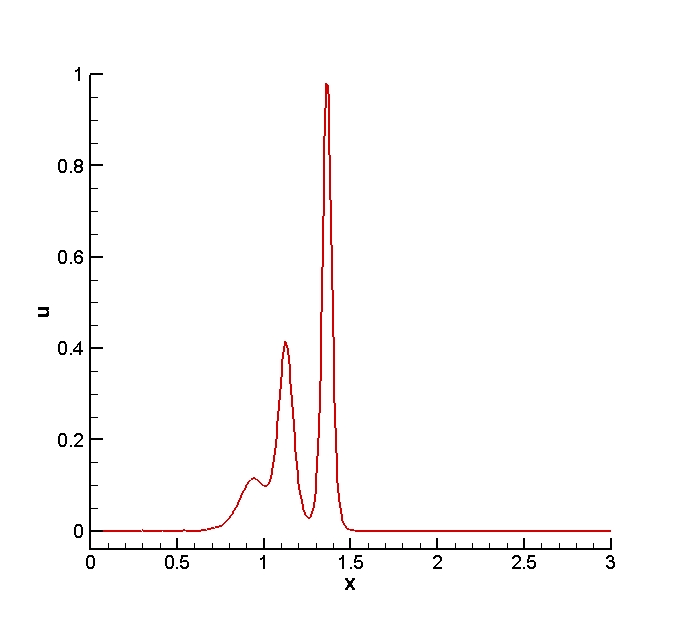}}\quad
    \subfigure[t=2]
    {\includegraphics[width=4.5cm]{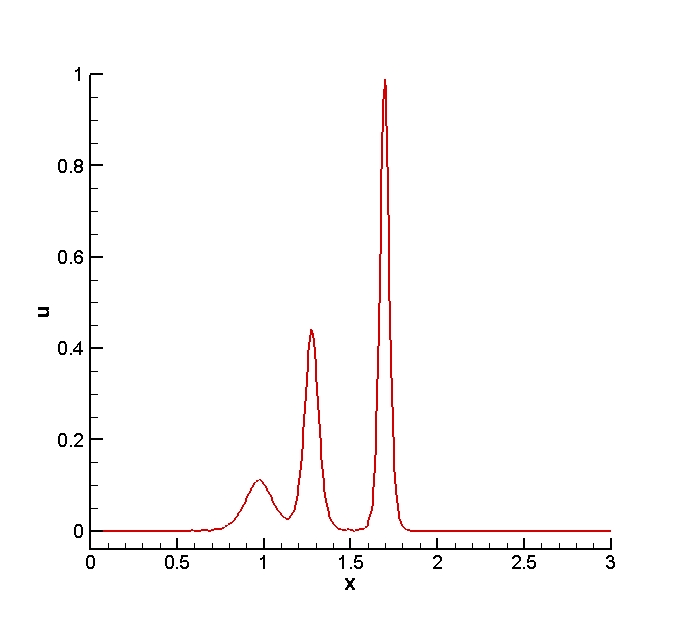}}
    }
    \caption{Example~\ref{exam4.3}: Triple soliton splitting case. Top: $P^2$ elements are used with 320 cells; Bottom: $P^4$ elements are used with 320 cells.}
    \label{fig:edge3}
    \end{center}
    \end{figure}

\begin{exam}{\em
\label{exam4.3}

In this example we compute several classical soliton solutions of the KdV equation
\[
u_t+(\frac{u^2}{2})_x+\epsilon u_{xxx}=0 .
\]
 Here, $f(u)=\frac{u^2}{2},\quad r(u)=u,\quad g(r)=\epsilon r$.

The single soliton case has the initial condition
\[
u(x,0)=3c\hbox {sech}^2(k(x-x_0))
\]
with $c=0.3$, $x_0=0.5$, $k=\frac{1}{2}\sqrt{\frac{c}{\epsilon}}$, and $\epsilon=5\times10^{-4}$.
The solution is computed in $x\in(0,2)$ with periodic boundary conditions and shown in Fig.~\ref{fig:edge1}.

The double soliton collision case has the initial condition
\[
u(x,0)=3c_1\hbox {sech}^2(k_1(x-x_1))+3c_2\hbox {sech}^2(k_2(x-y))
\]
with $c_1=0.3$, $c_2=0.1$, $x_1=0.4$, $y=0.8$, $k_i=\frac{1}{2}\sqrt{\frac{c_i}{\epsilon}}$ for $i=1,2$, and
$\epsilon=4.84\times10^{-4}$. The solution shown in Fig.~\ref{fig:edge2} is computed
in $x\in (0,2)$ with a periodic boundary condition.

The triple soliton splitting case has the initial condition
\[
u(x,0)=\frac{2}{3}\hbox {sech}^2(\frac{x-1}{\sqrt{108\epsilon}})
\]
with $\epsilon=10^{-4}$. The solution shown in Fig.~\ref{fig:edge3} is computed in $x\in(0,3)$ with
a periodic boundary condition.
}\end{exam}

\begin{figure}[hbtp]
  \begin{center}
  \mbox{\subfigure[$\epsilon=10^{-4}$, $t=0.5$, $n=200$]
  {\includegraphics[width=5.5cm]{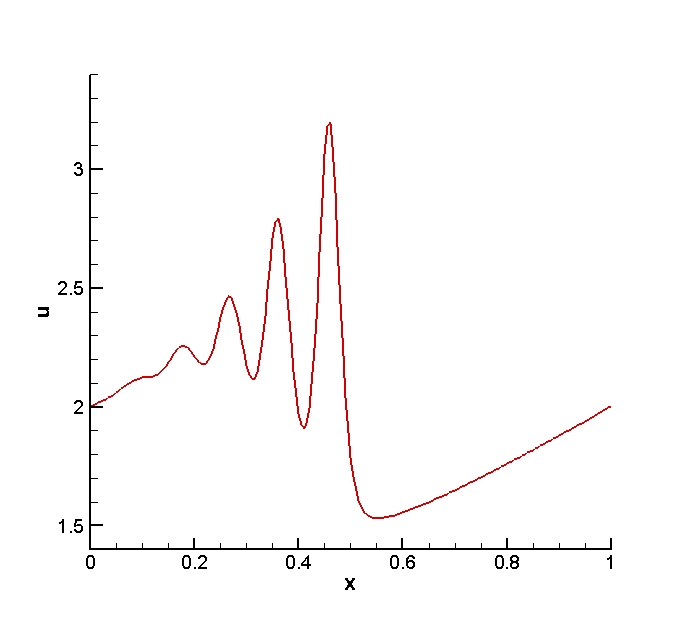}}\quad
    \subfigure[$\epsilon=10^{-5}$, $t=0.5$, $n=200$]
    {\includegraphics[width=5.5cm]{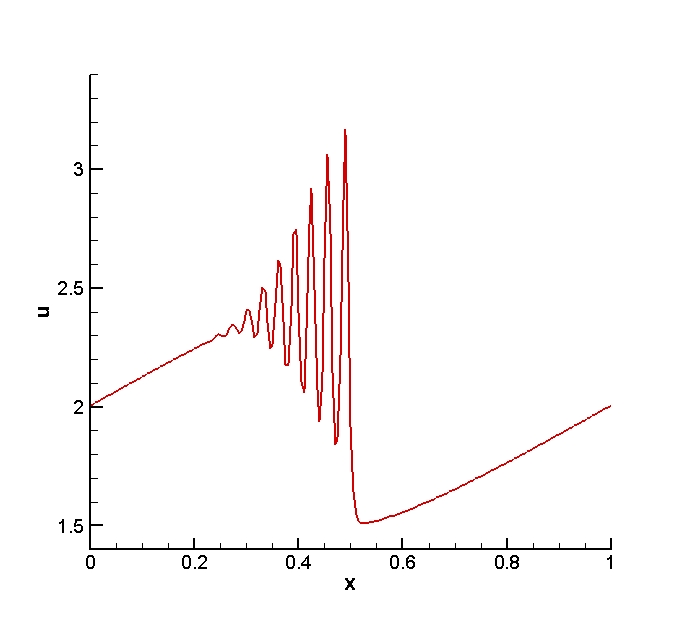}}
    }
    \mbox{\subfigure[$\epsilon=10^{-6}$, $t=0.5$, $n=1000$]
  {\includegraphics[width=5.5cm]{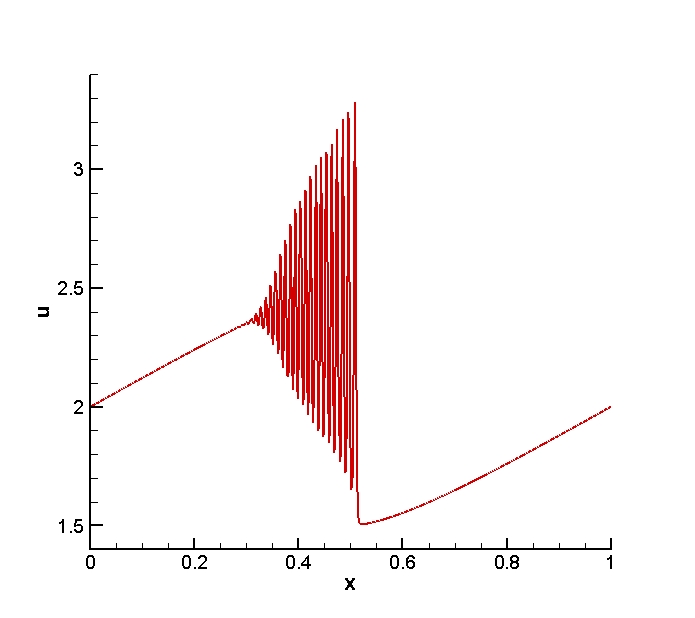}}\quad
    \subfigure[$\epsilon=10^{-7}$, $t=0.5$, $n=3400$]
    {\includegraphics[width=5.5cm]{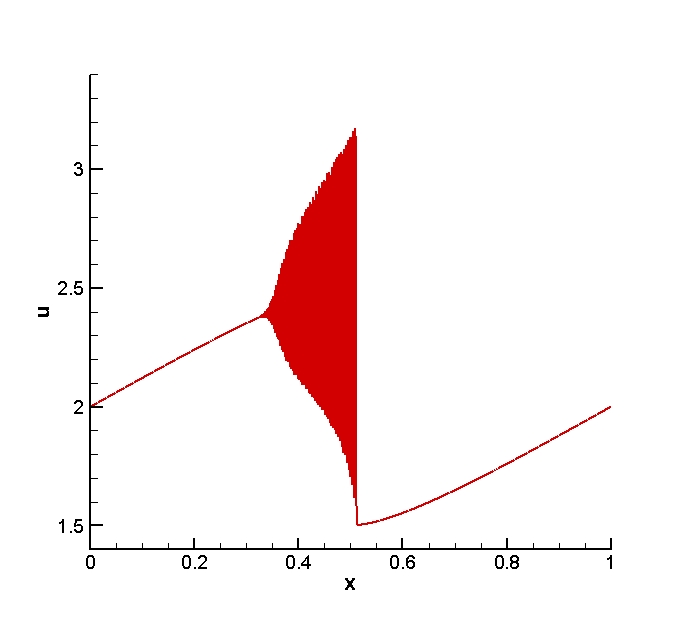}}
    }
    \caption{Example~\ref{exam4.4}: Zero dispersion limit of conservation laws. $P^2$ elements are used.}
    \label{fig:edge4}
    \end{center}
    \end{figure}

\begin{figure}[hbtp]
  \begin{center}
  \mbox{\subfigure[$\epsilon=10^{-4}$, $t=0.5$, $n=200$]
  {\includegraphics[width=5.5cm]{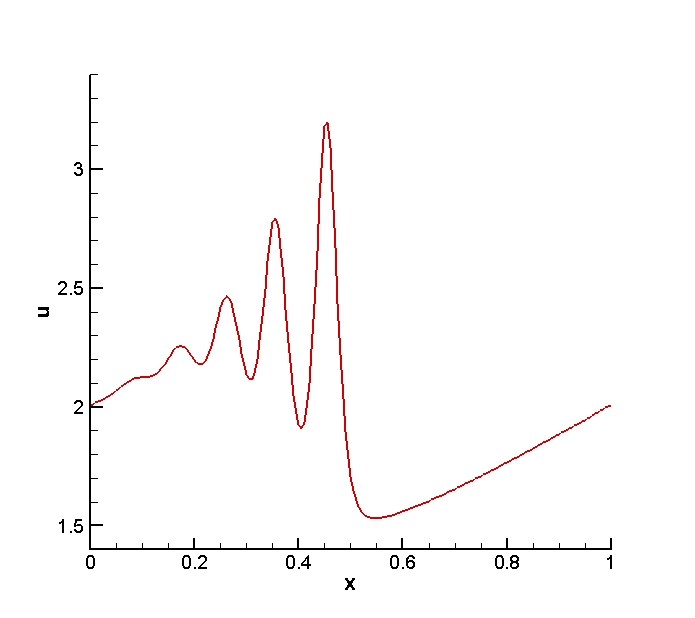}}\quad
    \subfigure[$\epsilon=10^{-5}$, $t=0.5$, $n=200$]
    {\includegraphics[width=5.5cm]{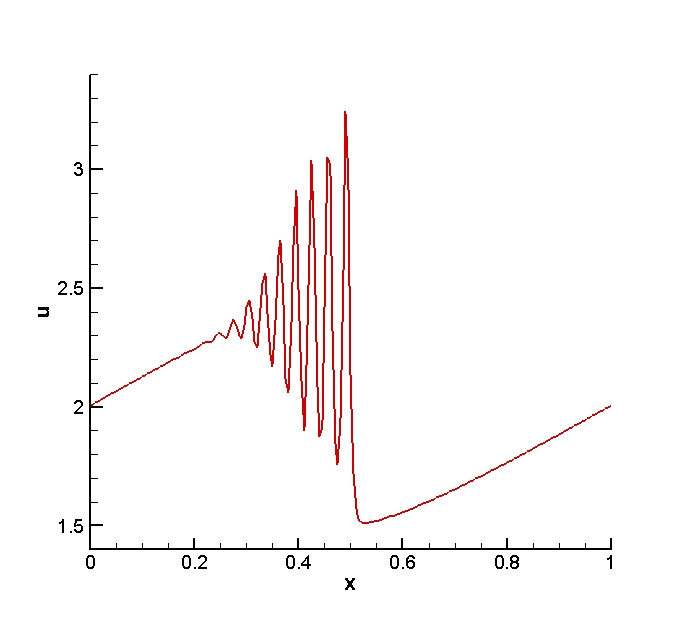}}
    }
    \mbox{\subfigure[$\epsilon=10^{-6}$, $t=0.5$, $n=1500$]
  {\includegraphics[width=5.5cm]{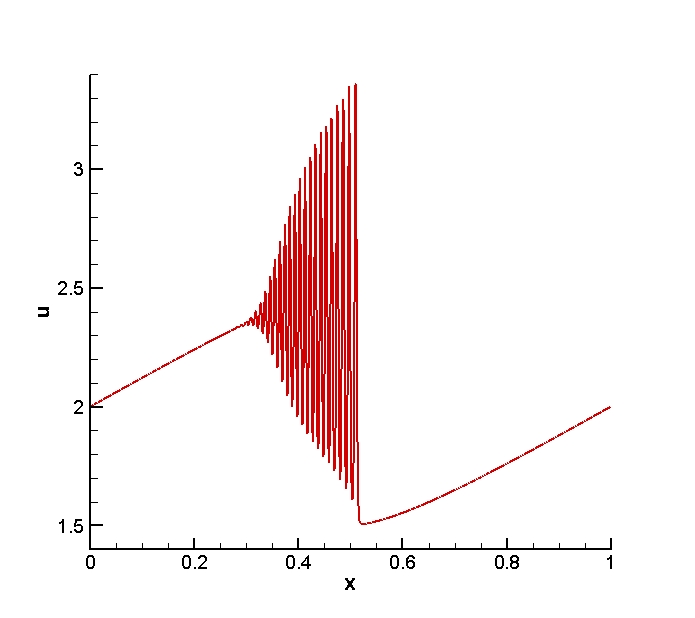}}\quad
    \subfigure[$\epsilon=10^{-7}$, $t=0.5$, $n=1800$]
    {\includegraphics[width=5.5cm]{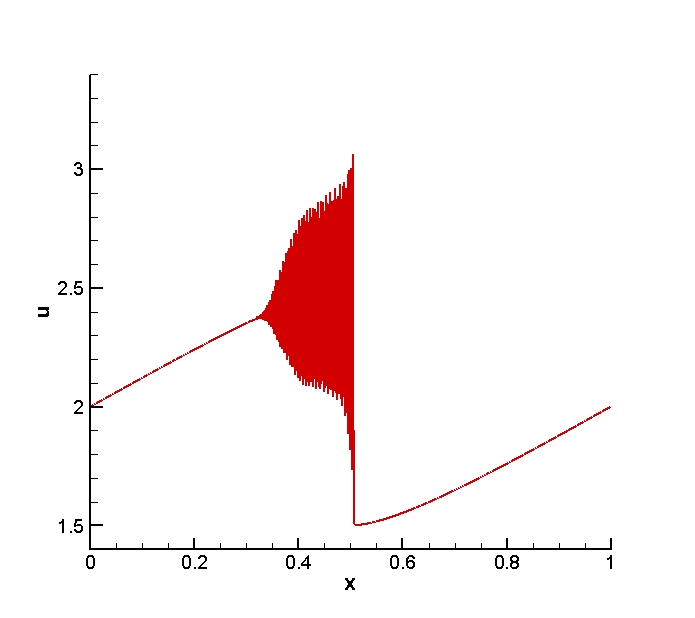}}
    }
    \caption{Example~\ref{exam4.4}: Zero dispersion limit of conservation laws. $P^4$ elements are used.}
    \label{fig:edge5}
    \end{center}
    \end{figure}

\begin{exam}{\em
\label{exam4.4}
We compute in this example the KdV zero dispersion limit of conservation laws. The equation is
\[
u_t+(\frac{u^2}{2})_x+\epsilon u_{xxx}=0,\quad x \in (0,1)
\]
subject to the periodic boundary condition and the initial condition
\[
u(x,0)=2+0.5\hbox{sin}(2\pi x).
\]
We compute the solution to $t=0.5$ with $\epsilon=10^{-4}$, $10^{-5}$, $10^{-6}$, and $10^{-7}$.
These numerical solutions are shown in Figs.~\ref{fig:edge4} and \ref{fig:edge5}.
}\end{exam}

\subsection{Two dimensional examples}

\begin{exam}{\em
\label{exam4.5}
In this example we compute the solution of the linear equation
\[
u_t+u_{xxx}+u_{yyy}=0, \quad (x,y) \in (0,2 \pi) \times (0, 2\pi)
\]
with the initial condition $u(x,y,0)=\sin(x+y)$ and the periodic boundary condition in both directions. Here, $f_1(u)=f_2(u)=0,\quad g_{12}=g_{21}=0,\quad g_{11}(r)=g_{22}(r)=r,\quad r_1(u)=r_2(u)=u.$
The exact solution is given by $u(x,y,t)=\sin(x+y-2t)$. The $L^1$, $L^2$, and $L^{\infty}$ norm of the error
and the convergence order of the scheme are shown in Table~\ref{ex3}. It can be seen clearly
that the scheme with $P^k$ elements in DG gives the $(k+1)^{\text{th}}$ order of convergence.
}\end{exam}

\begin{table}
\caption{Example~\ref{exam4.5}: periodic boundary conditions and $T=1$.}

\renewcommand{\multirowsetup}{\centering}
\begin{center}
\begin{tabular}{|c|c|c|c|c|c|c|}
\hline
$k$ & $n\times n$           & $10\times 10$ & $20\times 20$ & $30\times 30$ & $40\times 40$ & $50\times 50$\\ 
\hline
\multirow{6}{1cm}{2}

 & $L^1$       & 6.411e-3 & 4.351e-4 & 1.058e-4 & 4.123e-5  & 2.030e-5      \\
 &  Order      & \quad    & 3.88     & 3.49     & 3.28      &  3.18     \\
 &$L^2$       & 6.360e-3 & 4.563e-4 & 1.131e-4 & 4.454e-5   & 2.207e-6  \\
 &Order       & \quad    & 3.80     & 3.44     & 3.24       &  3.15     \\
 &$L_{\infty}$&7.974e-3  &5.973e-4  & 1.521e-4 & 6.093e-5  &  3.041e-5  \\
 &Order    & \quad       & 3.74    & 3.37     & 3.18      &   3.11     \\
\hline
\multirow{6}{1cm}{3}
 & $L^1$       & 5.993e-4 & 1.537e-5 & 1.839e-6 & 4.445e-7 & 1.571e-7 \\
 &  Order      & \quad    & 5.29     & 5.24     & 4.94      &  4.66\\
 &$L^2$       & 5.923e-4 & 1.572e-5 & 1.950e-6 & 4.821e-7 &1.726e-7  \\
 &Order       & \quad    & 5.24     & 5.15     & 4.86      & 4.60  \\
 &$L_{\infty}$&6.745e-4  &1.969e-5  & 2.703e-6 & 7.006e-7  & 2.554e-7 \\
 &Order    & \quad       & 5.10     & 4.90    & 4.69      &   4.52\\
\hline
\end{tabular}
\end{center}
\label{ex3}
\end{table}


\begin{exam}{\em
\label{exam4.6}
In this example we consider the Zakharov-Kuznetsov (ZK) equation
  \begin {align}
  u_t+uu_x+\varepsilon (u_{xxx}+u_{yyx})=0 .
  \label{ex6}
  \end{align}
  Here, $f_1(u)=\frac{u^2}{2},\quad f_2(u)=0,\quad g_{12}=g_{22}=0,\quad g_{11}(r)=g_{21}(r)=\epsilon r,\quad r_1(u)=r_2(u)=u.$
The steady progressive wave solution is of the form
  \begin {align}
  u(x,y,t)=3c\hbox{sech}^2(0.5\sqrt{{\frac{c}{\varepsilon}}}((x-ct-x_0)\hbox{cos}\theta+(y-y_0)\hbox{sin}\theta)),
  \end{align}
where $\theta$ is an inclined angle with respect to the $x$-axis and $(x_0,y_0)$ is the location of the peak of the initial $u$.
We can see in Tables~\ref{ex4} and \ref{ex5} that the method with $P^k$ elements gives the $(k+1)^{\text{th}}$
order of convergence. The computational domains are $(-16,16)\times (-16,16)$ in Table \ref{ex4}
and $(0,32)\times(0,16)$ in Table \ref{ex5}, respectively.
}\end{exam}

\begin{table}
\caption{Example~\ref{exam4.6}: periodic boundary condition in both directions, $c=0.01$,
$\varepsilon=0.01$, $\theta =0$, $x_0=0$, $y_0=0$, $T=1$.}

\renewcommand{\multirowsetup}{\centering}
\begin{center}
\begin{tabular}{|c|c|c|c|c|c|c|}
\hline
$k$ & $n\times n$           & $40\times 40$ & $50\times 50$ & $60\times 60$ & $70\times 70$ & $80\times 80$\\ 
\hline
\multirow{6}{1cm}{2}

 & $L^1$       & 1.515e-6 & 4.671e-7 & 2.202e-7 & 1.120e-7  & 6.404e-8      \\
 &  Order      & \quad    & 5.27     & 4.12     & 4.39      &  4.19     \\
 &$L^2$       & 4.128e-6 & 1.527e-6 & 7.131e-7 & 3.887e-7   & 2.189e-7  \\
 &Order       & \quad    & 4.46     & 4.18    & 3.94       &  4.30     \\
 &$L_{\infty}$&1.747e-5  &7.278e-6  & 4.573e-6 & 2.659e-6  &  1.449e-6  \\
 &Order    & \quad       & 3.92    & 2.55     & 3.52      &   4.55   \\
\hline
\multirow{6}{1cm}{3}
 & $L^1$       & 1.046e-6 & 2.796e-7 & 1.249e-7 & 5.814e-8 & 2.785e-8 \\
 &  Order      & \quad    & 5.91     & 4.42     & 4.96      &  5.51\\
 &$L^2$       & 2.376e-6 & 8.512e-7 & 3.663e-7 & 1.869e-7 &9.563e-8  \\
 &Order       & \quad    & 4.60    & 4.62     & 4.37      & 5.02  \\
 &$L_{\infty}$&8.766e-6  &4.474e-6  & 1.866e-6 & 1.112e-6  & 5.968e-7 \\
 &Order    & \quad       & 3.01     & 4.80     & 3.36      &   4.66\\
\hline
\end{tabular}
\end{center}
\label{ex4}
\end{table}

\begin{table}
\caption{Example~\ref{exam4.6}: periodic boundary condition in the $x$-direction and Dirichlet boundary condition
in the $y$-directions, $c=0.01$, $\varepsilon=0.01$, $\theta =\frac{\pi}{12}$, $x_0=16$, $y_0=8$, and $T=1$.}

\renewcommand{\multirowsetup}{\centering}
\begin{center}
\begin{tabular}{|c|c|c|c|c|c|c|}
\hline
$k$ & $n\times n$           & $40\times 40$ & $50\times 50$ & $60\times 60$ & $70\times 70$ & $80\times 80$\\ 
\hline
\multirow{6}{1cm}{2}

 & $L^1$       & 1.583e-6 & 5.111e-7 & 2.366e-7 & 1.088e-7  & 5.919e-8      \\
 &  Order      & \quad    & 5.07     & 4.22     & 5.04      &  4.56     \\
 &$L^2$       & 5.009e-6 & 1.802e-6 & 8.382e-7 & 3.922e-7   & 1.966e-7  \\
 &Order       & \quad    & 4.58     & 4.20     & 4.93       &  5.17     \\
 &$L_{\infty}$&3.623e-5  &1.538e-5  & 6.424e-6 & 3.400e-6  &  1.592e-6  \\
 &Order    & \quad       & 3.84    & 4.79   & 4.13      &   5.68     \\
\hline
\multirow{6}{1cm}{3}
 & $L^1$       & 1.523e-6 & 3.564e-7 & 1.800e-7 & 6.543e-8 & 2.846e-8 \\
 &  Order      & \quad    & 6.51    & 3.75     & 6.56     &  6.23\\
 &$L^2$       & 4.014e-6 & 1.164e-6 & 5.582e-7 & 2.337e-7 &1.009e-7  \\
 &Order       & \quad    & 5.55     & 4.03     & 5.65      & 6.29  \\
 &$L_{\infty}$&3.094e-5  &9.814e-6  & 4.009e-6 & 2.211e-6  & 8.719e-7 \\
 &Order    & \quad       & 5.15     & 4.91     & 3.86      &   6.97\\
\hline
\end{tabular}
\end{center}
\label{ex5}
\end{table}

\begin{figure}[hbtp]
  \begin{center}
  \mbox{\subfigure[$t=0$]
  {\includegraphics[width=5.5cm]{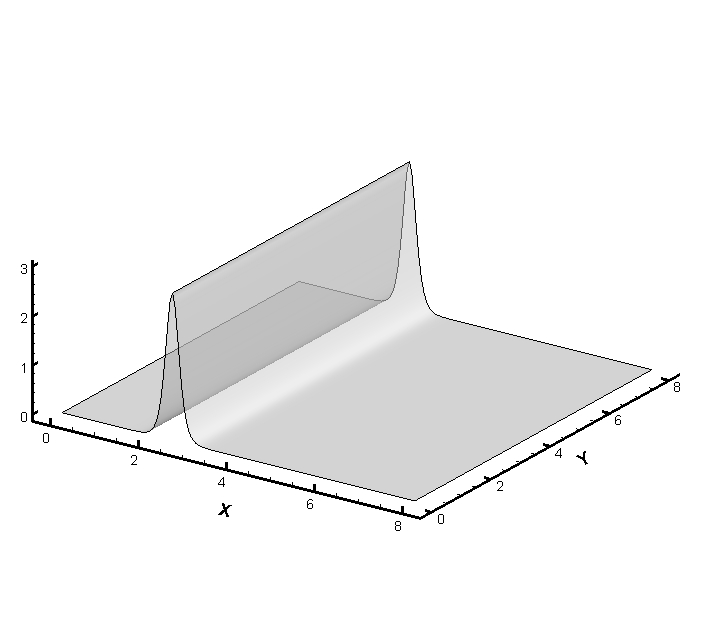}}\quad
    \subfigure[$t=1$]
    {\includegraphics[width=5.5cm]{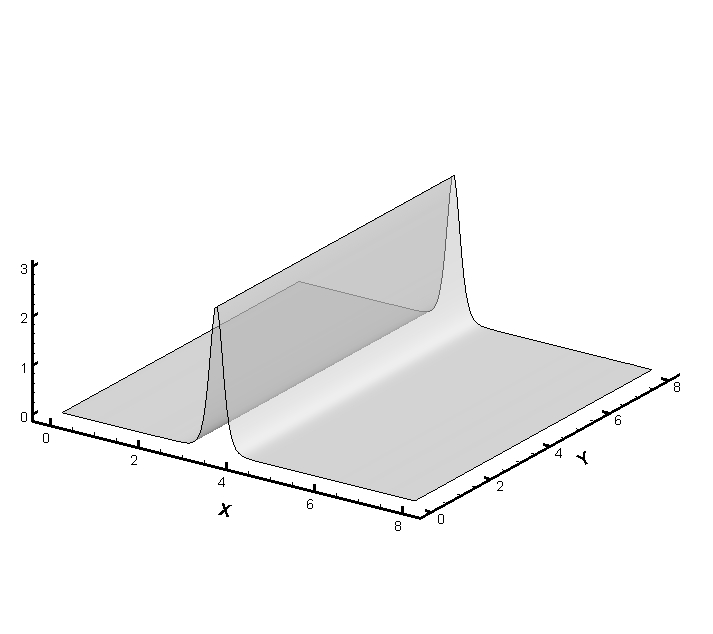}}
    }
    \mbox{\subfigure[$t=2$]
  {\includegraphics[width=5.5cm]{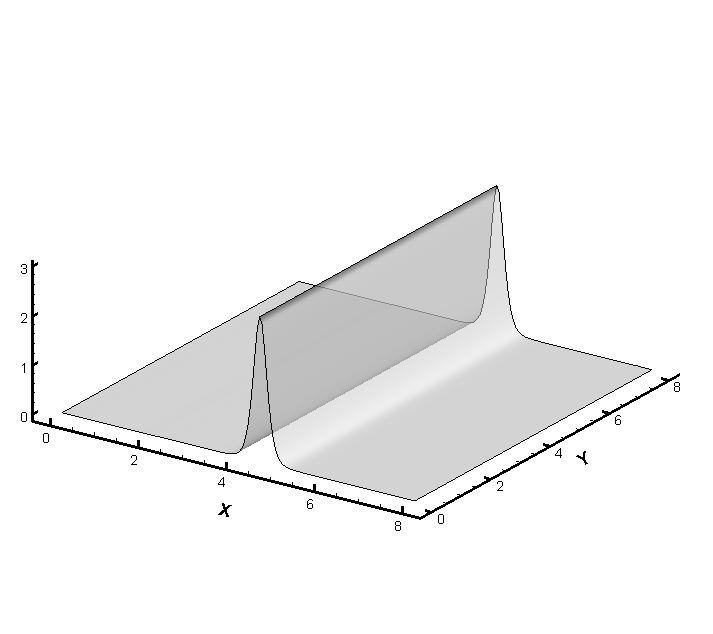}}\quad
    \subfigure[$t=3$]
    {\includegraphics[width=5.5cm]{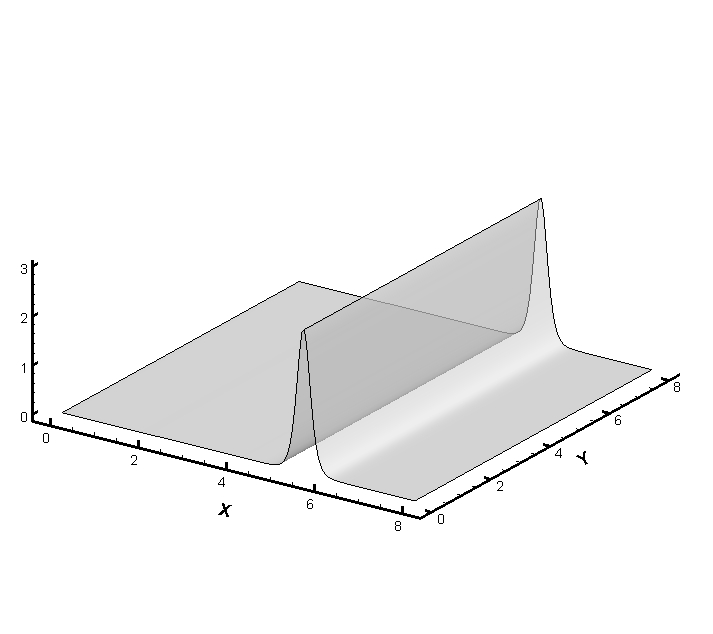}}
    }
    \caption{Example~\ref{exam4.7}: The single steady progressive wave propagation with the initial condition (\ref{inti0}). $c=1$, $\varepsilon=0.01$, $\theta=0$, $x_0=2.5$, $y_0=4$. Periodic boundary condition in both coordinate
    directions in $(0,8)\times(0,8)$. $P^2$ elements with $150\times150$ uniform cells are used.}
    \label{fig:edge6}
    \end{center}
    \end{figure}

\begin{figure}[hbtp]
  \begin{center}
  \mbox{\subfigure[$t=0$]
  {\includegraphics[width=5.5cm]{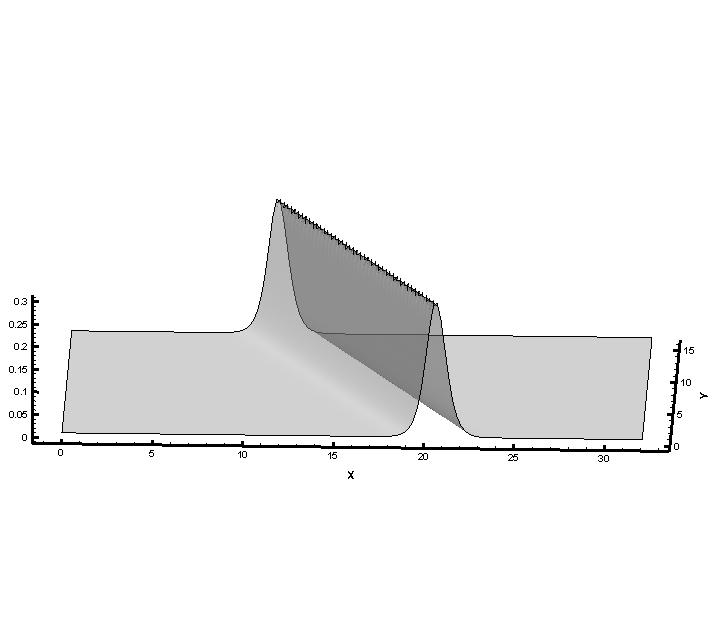}}\quad
    \subfigure[$t=15$]
    {\includegraphics[width=5.5cm]{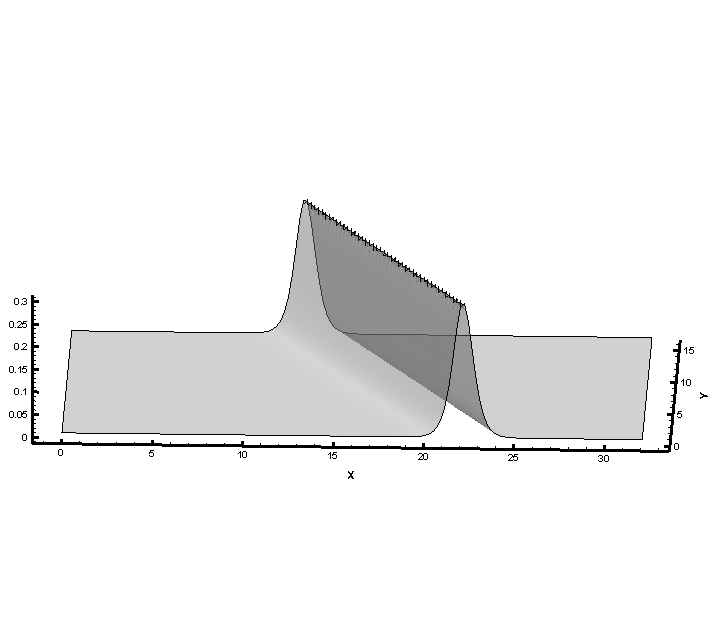}}
    }
    \mbox{\subfigure[$t=30$]
  {\includegraphics[width=5.5cm]{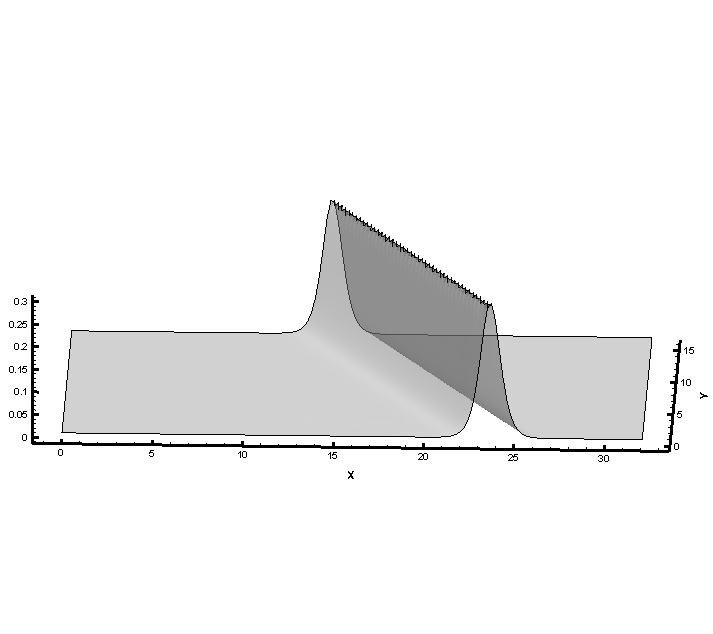}}\quad
    \subfigure[$t=45$]
    {\includegraphics[width=5.5cm]{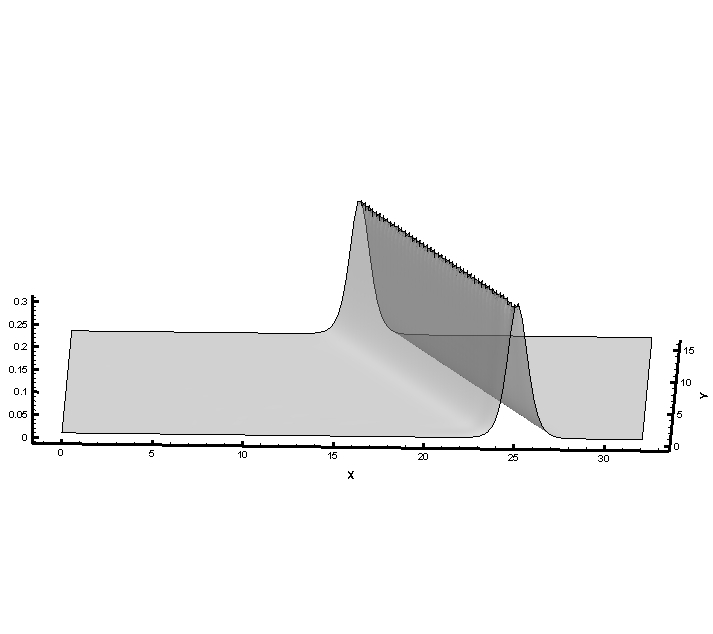}}
    }
    \caption{Example~\ref{exam4.7}: The single steady progressive wave propagation with the initial condition (\ref{inti0}). $c=0.1$, $\varepsilon=0.01$, $\theta=\frac{\pi}{6}$, $x_0=16$, $y_0=8$. Periodic boundary condition in the $x$-direction and Dirichlet boundary condition in the $y$-direction in $(0,32)\times(0,16)$. $P^2$ elements with $150\times150$ uniform cells are used.}
    \label{fig:edge7}
    \end{center}
    \end{figure}

\begin{figure}[hbtp]
  \begin{center}
  \mbox{\subfigure[$t=0$]
  {\includegraphics[width=5.5cm]{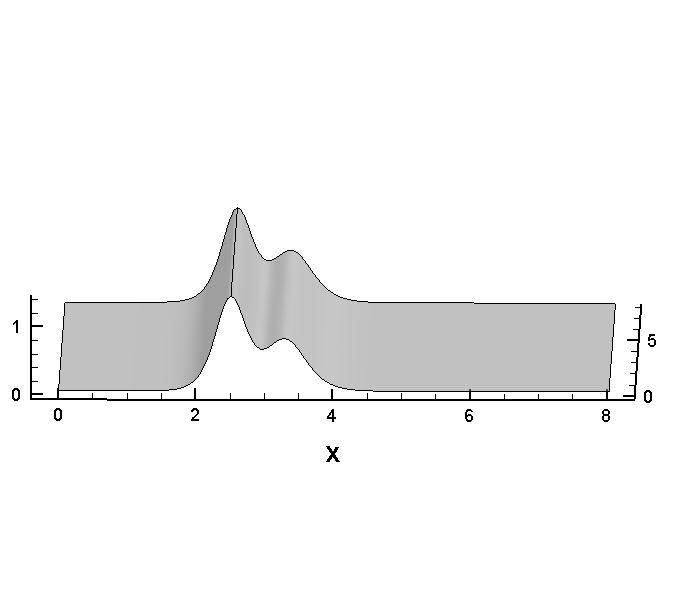}}\quad
    \subfigure[$t=1$]
    {\includegraphics[width=5.5cm]{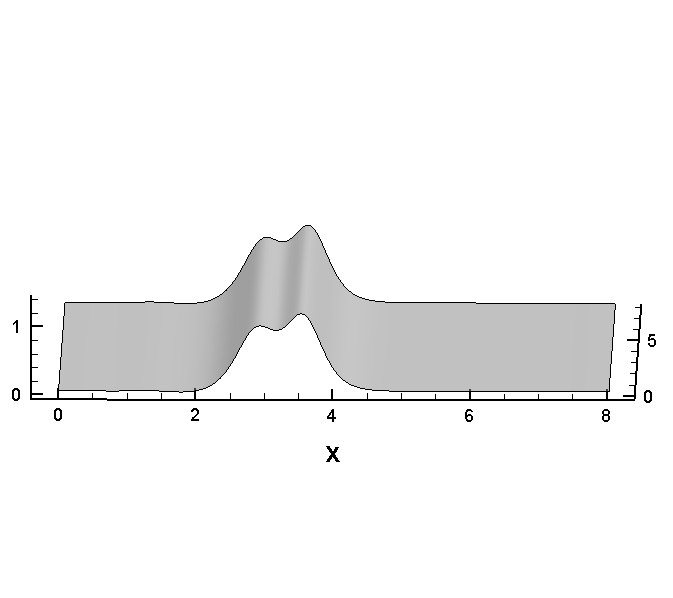}}
    }
    \mbox{\subfigure[$t=2$]
  {\includegraphics[width=5.5cm]{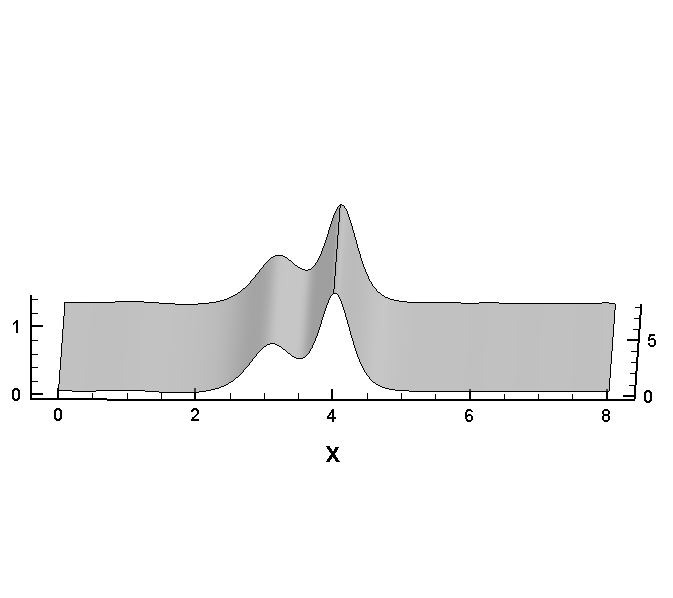}}\quad
    \subfigure[$t=3$]
    {\includegraphics[width=5.5cm]{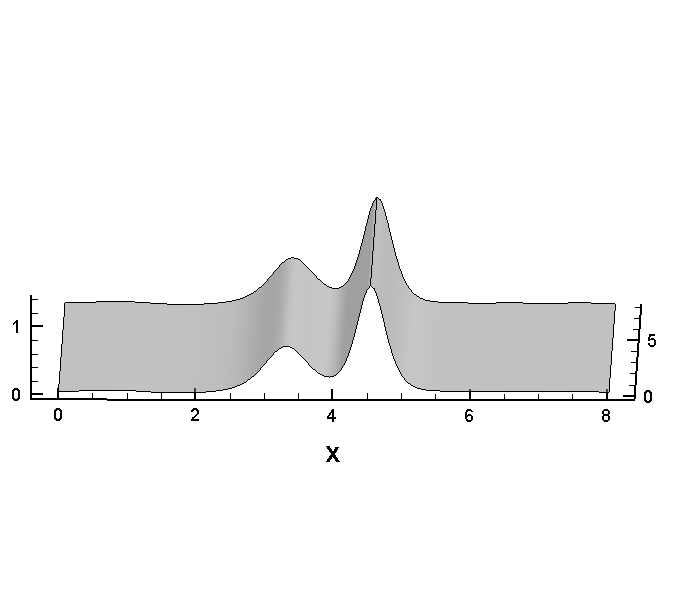}}
    }
    \caption{Example~\ref{exam4.7}: Double soliton collision profiles and periodic boundary condition in both  coordinate directions in $(0,8)\times(0,8)$. $P^2$ elements with $150\times150$ uniform cells are used.}
    \label{fig:edge8}
    \end{center}
    \end{figure}

\begin{exam}{\em
\label{exam4.7}
In this example we show the steady progressive wave propagation of the ZK equation (\ref{ex6}). First, we show the single steady progressive wave in Figs. \ref{fig:edge6} and \ref{fig:edge7} with the initial condition
\begin {align}
  u(x,y,0)=3c\hbox{sech}^2(0.5\sqrt{{\frac{c}{\varepsilon}}}((x-x_0)\hbox{cos}\theta+(y-y_0)\hbox{sin}\theta)) ,
  \label{inti0}
  \end{align}
  where $\theta$ is the inclined angle with respect to the $x$-axis. The periodic boundary condition in both $x$-
  and $y$-directions are used when $\theta =0$. Since the solution cannot be periodic in $y$-direction when
  $\theta\neq 0$, we use for this case the Dirichlet boundary condition in the $y$-direction
  and a periodic boundary condition in the $x$-direction.

  The double soliton collision case has the initial condition
  \begin {align}
  u(x,y,0)=\sum\limits_{j=1}^23c_j\hbox{sech}^2(0.5\sqrt{{\frac{c_j}{\varepsilon}}}((x-x_j)\hbox{cos}\theta+(y-y_j)\hbox{sin}\theta)),
  \end{align}
  where $c_1=0.45$, $c_2=0.25$, $\varepsilon=0.01$, $\theta=0$, $x_1=2.5$, $y_1=0$, $y=3.3$, $y_2=0$.
  The results with periodic boundary conditions in both coordinate directions in $(0,8)\times(0,8)$
  using $P^2$ elements with $150\times150$ uniform cells are shown in Fig.~\ref{fig:edge8}.

}\end{exam}


\begin{figure}[hbtp]
  \begin{center}
  \mbox{\subfigure[$t=0$]
  {\includegraphics[width=5.5cm]{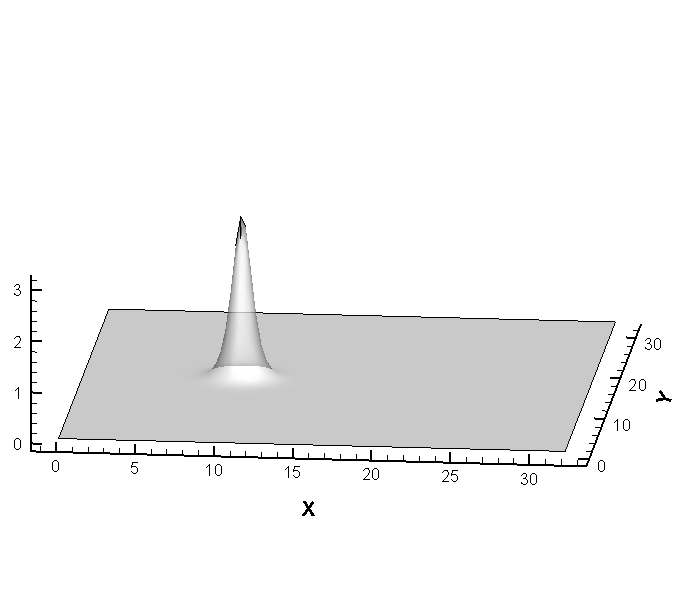}}\quad
    \subfigure[$t=1$]
    {\includegraphics[width=5.5cm]{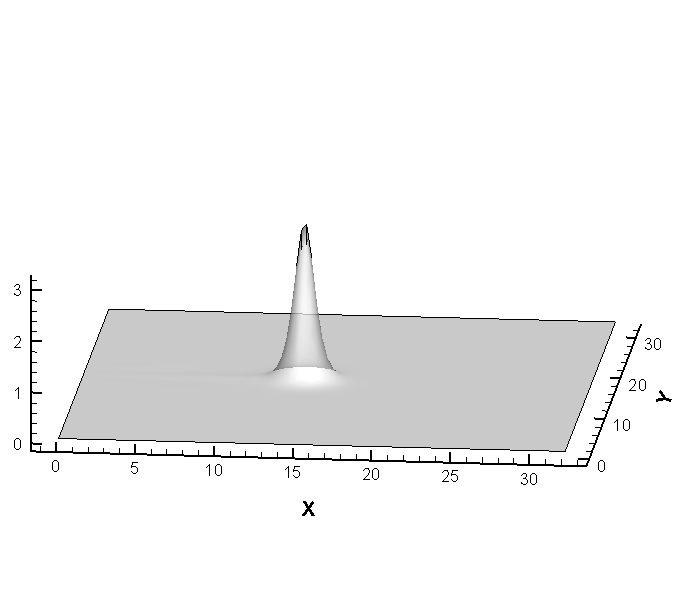}}
    }
    \mbox{\subfigure[$t=2$]
  {\includegraphics[width=5.5cm]{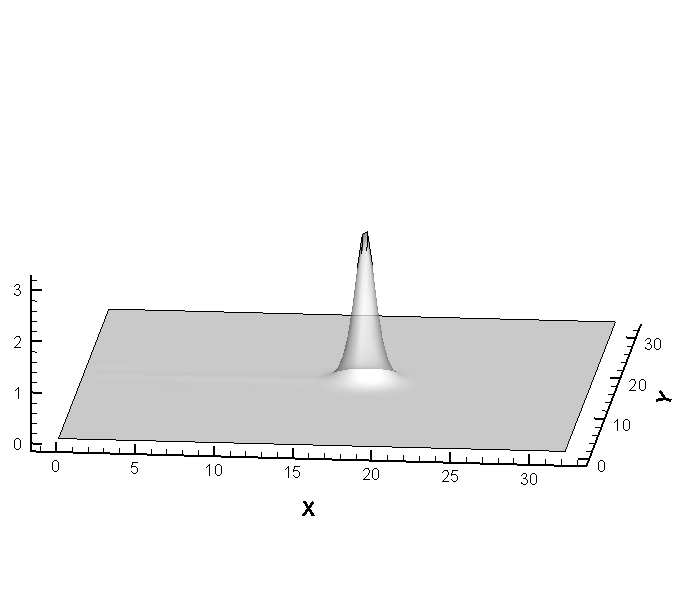}}\quad
    \subfigure[$t=3$]
    {\includegraphics[width=5.5cm]{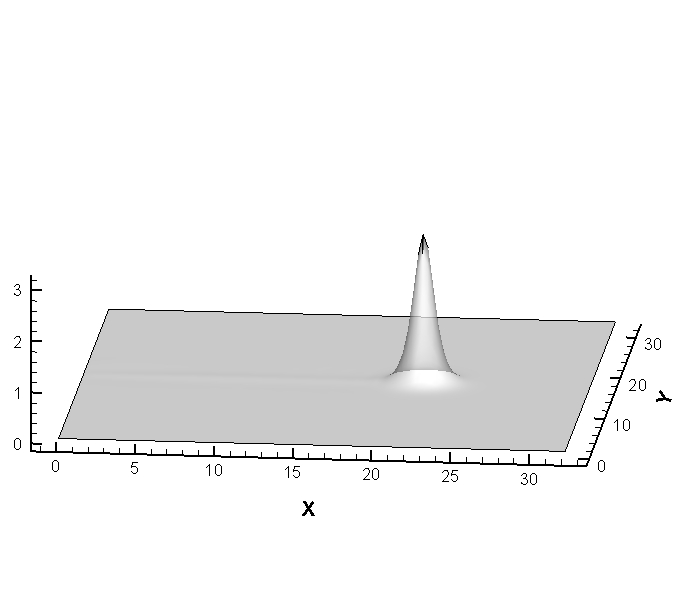}}
    }
    \caption{Example~\ref{exam4.8}: Evolution of a single bell-shaped pulse solution for (\ref{ex4.8}) with $c=4$,
    $x_0=10$, and $y_0=16$. The used domain is $(0,32)\times(0,32)$.
    $P^2$ elements with $100\times100$ uniform cells are used.}
    \label{fig:edge9}
    \end{center}
    \end{figure}

\begin{exam}{\em
\label{exam4.8}
In this example we show the numerical results for the equation
\begin {align}
  u_t+(3u^2)_x+u_{xxx}+u_{xyy}=0 .
  \label{ex4.8}
  \end{align}
Here, $f_1(u)=3u^2,\quad f_2(u)=0,\quad g_{12}=g_{22}=0,\quad g_{11}(u)=g_{21}(u)=u,\quad r_1(u)=r_2(u)=u.$ A  cylindrically symmetric solitary solution was obtained and its evolutions as well as interactions were investigated numerically in \cite{ER11}. This type of solitary solution, also called the bell-shaped pulse, can be expressed as
\begin {align}
  u(x,y,t)=\frac{c}{3}\sum\limits_{n=1}^{10}a_{2n}(\hbox{cos}(2n\, \hbox{arccot}(\frac{\sqrt{c}}{2}r))-1) ,
  \label{ex4.8.1}
  \end{align}
where $c$ is the velocity of the solitary wave solution and $r=\sqrt{(x-ct)^2+y^2}$.
The coefficients are collected in Table~\ref{ex7} . We use the Dirichlet boundary condition given by
the exact solution. The stable propagation of a single pulse is shown in Fig.~\ref{fig:edge9}.

Then, we proceed to show the collision of two pulses with the initial condition
\begin {align}
  u(x,y,0)=\sum\limits_{j=1}^2\frac{c_j}{3}\sum\limits_{n=1}^{10}a_{2n}(\hbox{cos}(2n\,\hbox{arccot}(\frac{\sqrt{c_j}}{2}r_j))-1),
  \label{ex4.8.2}
  \end{align}
  where $c_1$ and $c_2$ are the velocities of the solitary wave solutions,
 $r_i=\sqrt{(x-x_i)^2+(y-y_i)^2}$ $(i=1,2)$, and $(x_i,y_i)$'s are the locations of the peaks of $u$.
 When the centers of the two pulses are situated on the same line with $y=\text{const}$,
 the collision is  called a direct collision and otherwise called a deviated collision \cite{ER11}.
The numerical solutions obtained for the cases of a direct collision and a deviated
 collision of two pulses are shown in Figs.~\ref{fig:edge10} and \ref{fig:edge11}, respectively.

}\end{exam}

\begin{table}
\caption{Example~\ref{exam4.8}: Coefficients for the solitary wave solution of the equation (\ref{ex4.8}).}
\renewcommand{\multirowsetup}{\centering}
\begin{center}
\begin{tabular}{|c|c|c|c|c|c|c|}
\hline
 $n$  & $a_{2n}$ & $n$ & $a_{2n}$ \\
\hline
 1& -1.25529873  & 6 & -0.00281281  \\
 2&  0.21722635  & 7 & -0.00138352  \\
 3&  0.06452543  & 8 & -0.00070289  \\
 4&  0.00540862  & 9 & -0.00020451  \\
 5& -0.00332515  & 10& -0.00003053  \\
\hline
\end{tabular}
\end{center}
\label{ex7}
\end{table}

    \begin{figure}[hbtp]
  \begin{center}
  \mbox{\subfigure[$t=0$]
  {\includegraphics[width=5.5cm]{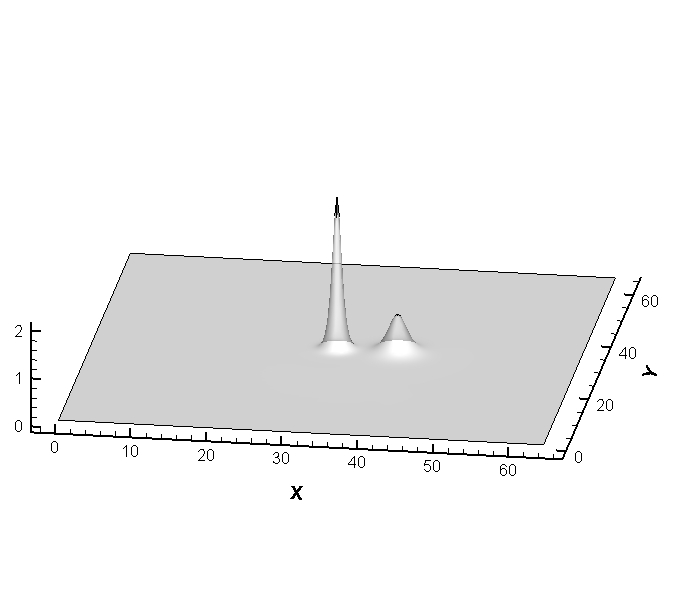}}\quad
    \subfigure[$t=0$]
    {\includegraphics[width=5.5cm]{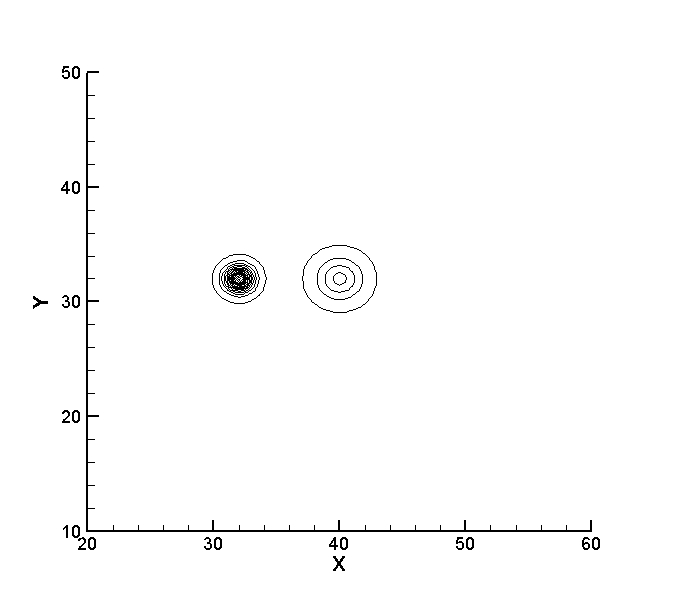}}
    }
    \mbox{\subfigure[$t=2$]
  {\includegraphics[width=5.5cm]{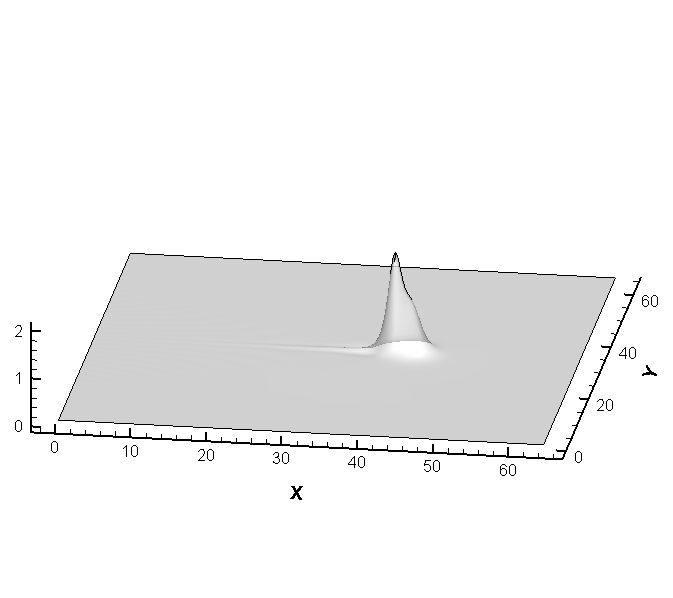}}\quad
    \subfigure[$t=2$]
    {\includegraphics[width=5.5cm]{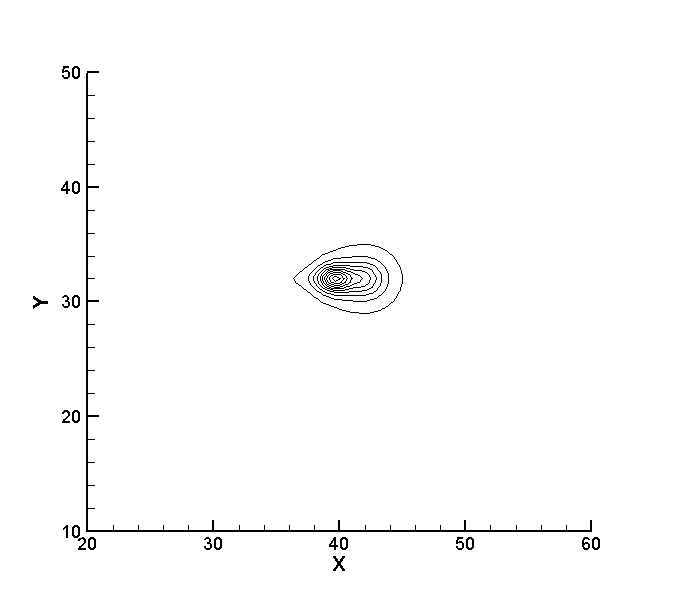}}
    }
    
    \caption{Example~\ref{exam4.8}: Direct collision of two bell-shaped pulses solution for (\ref{ex4.8}) with the initial condition (\ref{ex4.8.2}) and $c_1=4$, $c_2=1$, $x_1=32$, $y_1=32$, $y=40$, $y_2=32$, and
Dirichlet boundary condition. The computational domain is $(0,64)\times(0,64)$.
$P^2$ elements with $200\times200$ uniform cells are used.}
    \label{fig:edge10}
    \end{center}
    \end{figure}

    \addtocounter{figure}{-1}
    \begin{figure}[hbtp]
  \begin{center}
  \addtocounter{subfigure}{1}      

  \addtocounter{subfigure}{-1}     
  \addtocounter{subfigure}{4}
    \mbox{\subfigure[$t=4$]
  {\includegraphics[width=5.5cm]{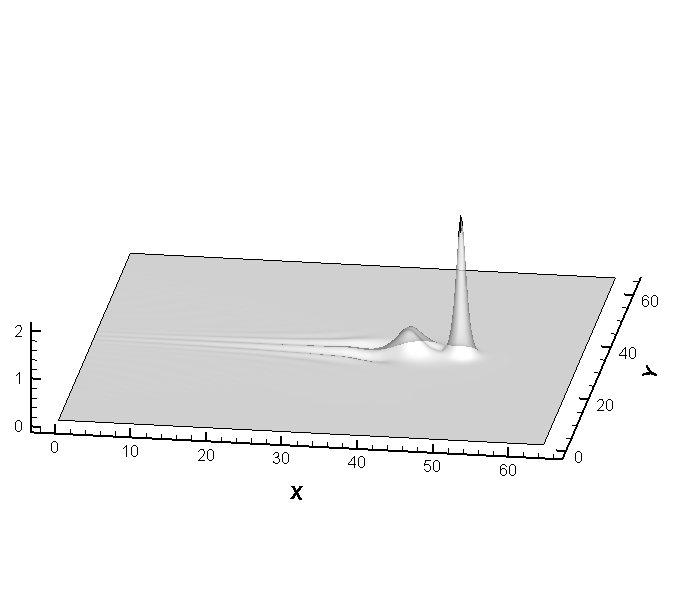}}\quad
    \subfigure[$t=4$]
    {\includegraphics[width=5.5cm]{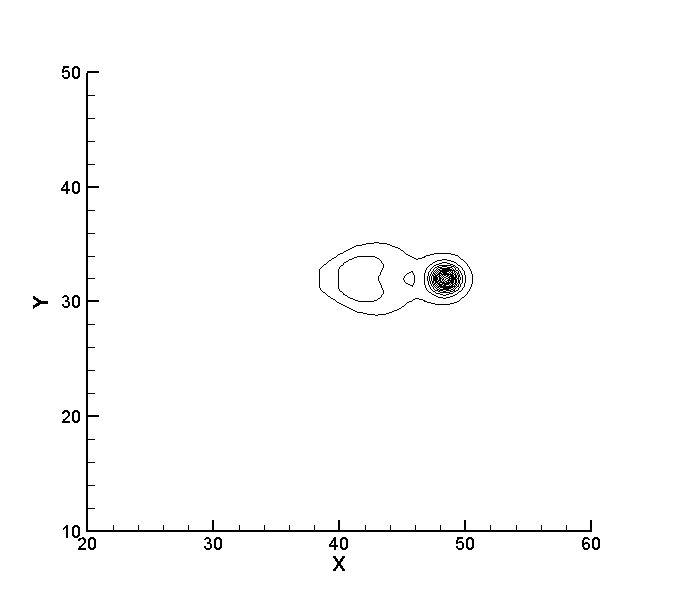}}
    }
     \mbox{\subfigure[$t=6$]
  {\includegraphics[width=5.5cm]{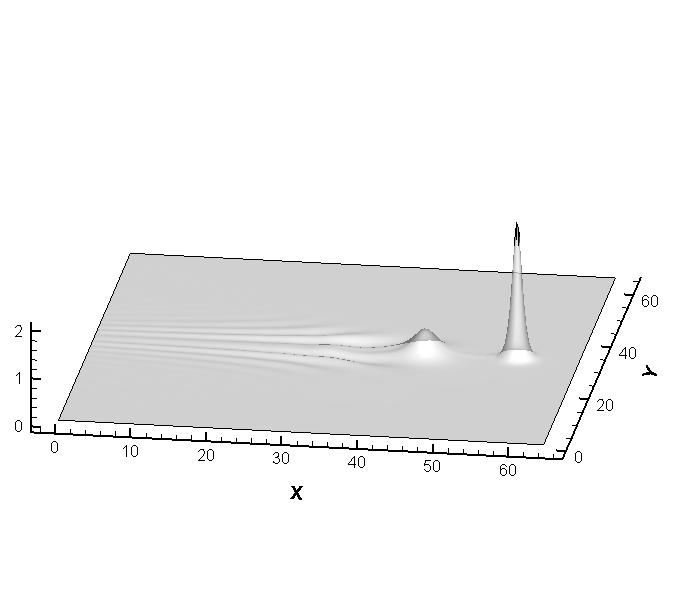}}\quad
    \subfigure[$t=6$]
    {\includegraphics[width=5.5cm]{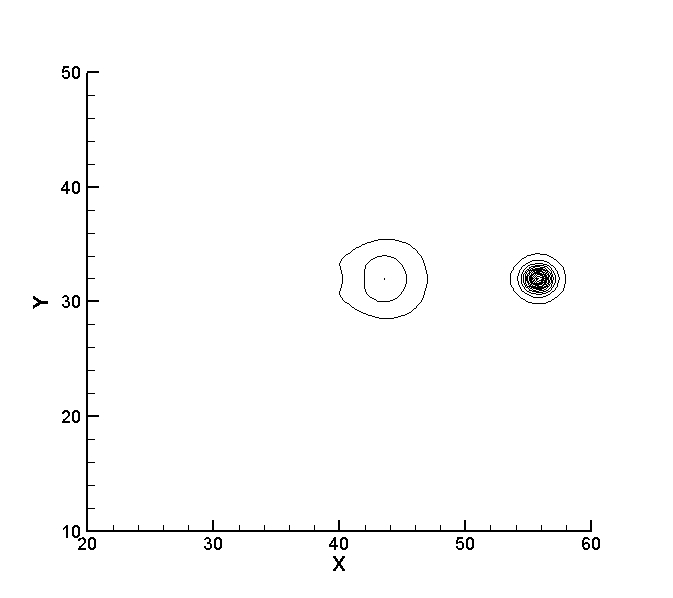}}
    }
    \caption{(continued)}
    \label{fig:edge10}
    \end{center}
    \end{figure}

\begin{figure}[hbtp]
  \begin{center}
  \mbox{\subfigure[$t=0$]
  {\includegraphics[width=5.5cm]{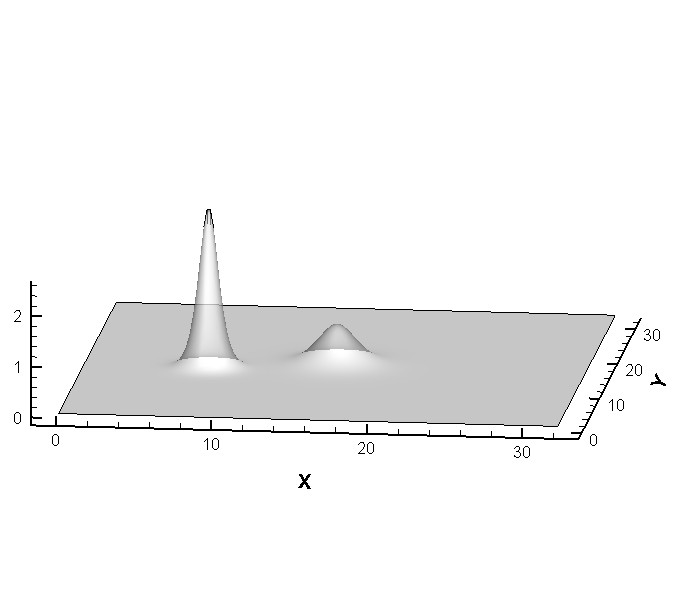}}\quad
    \subfigure[$t=0$]
    {\includegraphics[width=5.5cm]{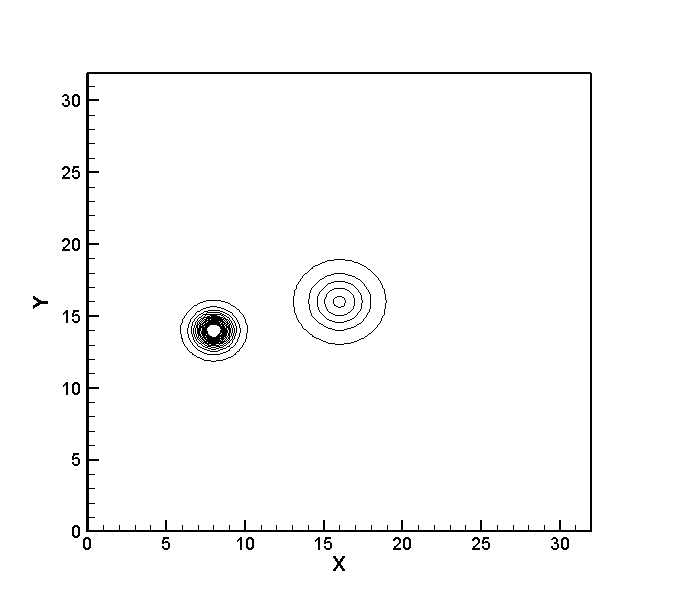}}
    }
    \mbox{\subfigure[$t=2$]
  {\includegraphics[width=5.5cm]{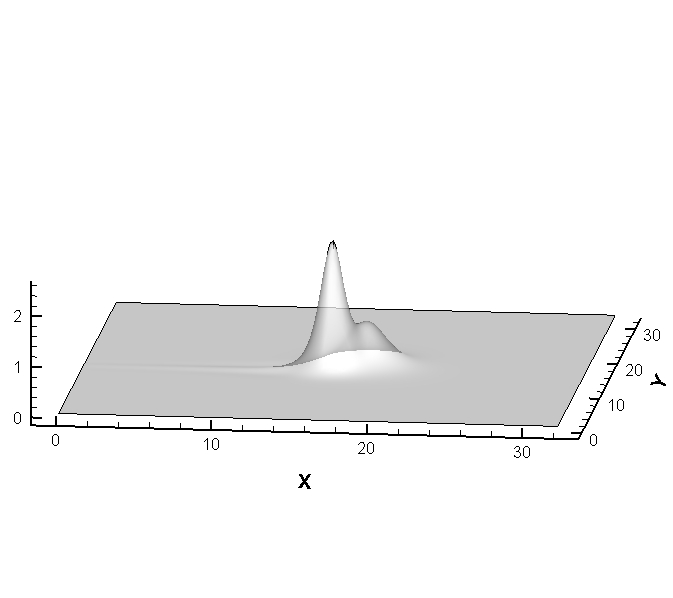}}\quad
    \subfigure[$t=2$]
    {\includegraphics[width=5.5cm]{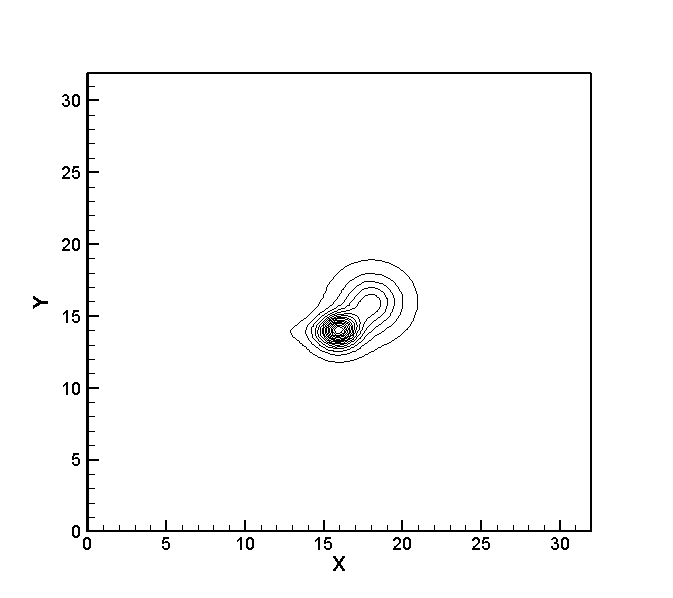}}
    }

    \caption{Example~\ref{exam4.8}: Deviated collision of two bell-shaped pulses solution for (\ref{ex4.8}) with the initial condition (\ref{ex4.8.2}) and $c_1=4$, $c_2=1$, $x_1=8$, $y_1=14$, $y=16$, $y_2=16$. and
    Dirichlet boundary condition. The computational domain is $(0,32)\times(0,32)$.
    $P^2$ elements with $150\times150$ uniform cells are used.}

    \label{fig:edge11}
    \end{center}
    \end{figure}
\addtocounter{figure}{-1}
    \begin{figure}[hbtp]
  \begin{center}
  \addtocounter{subfigure}{1}      

  \addtocounter{subfigure}{-1}     
  \addtocounter{subfigure}{4}
  \centering

  \mbox{\subfigure[$t=3$]
  {\includegraphics[width=5.5cm]{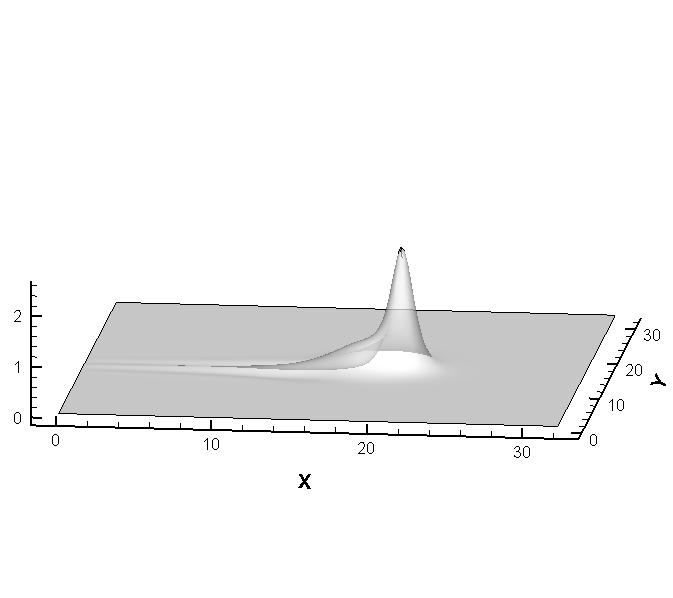}}\quad
    \subfigure[$t=3$]
    {\includegraphics[width=5.5cm]{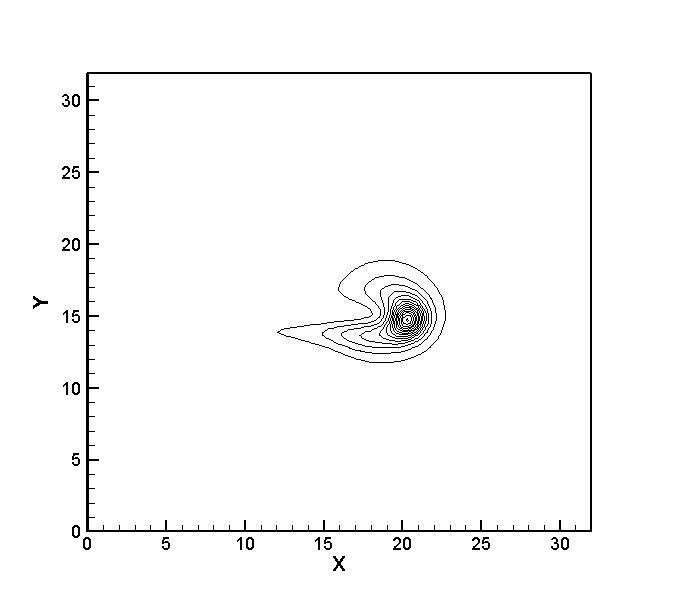}}
    }
  \mbox{\subfigure[$t=5$]
  {\includegraphics[width=5.5cm]{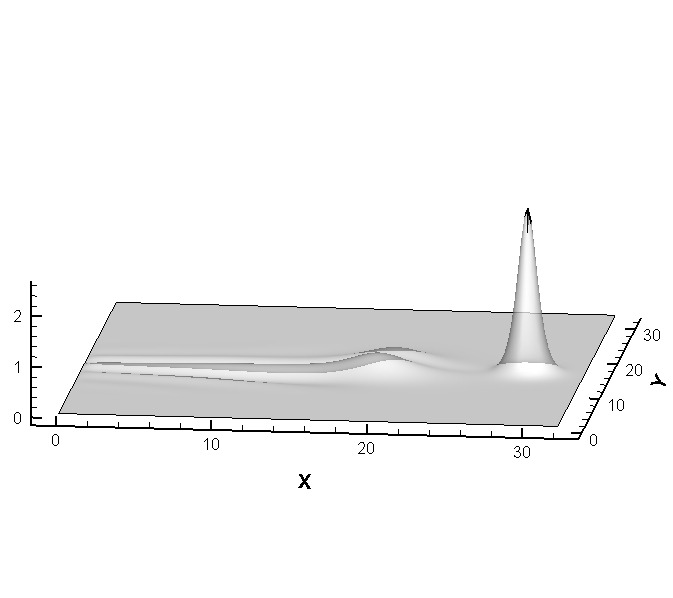}}\quad
    \subfigure[$t=5$]
    {\includegraphics[width=5.5cm]{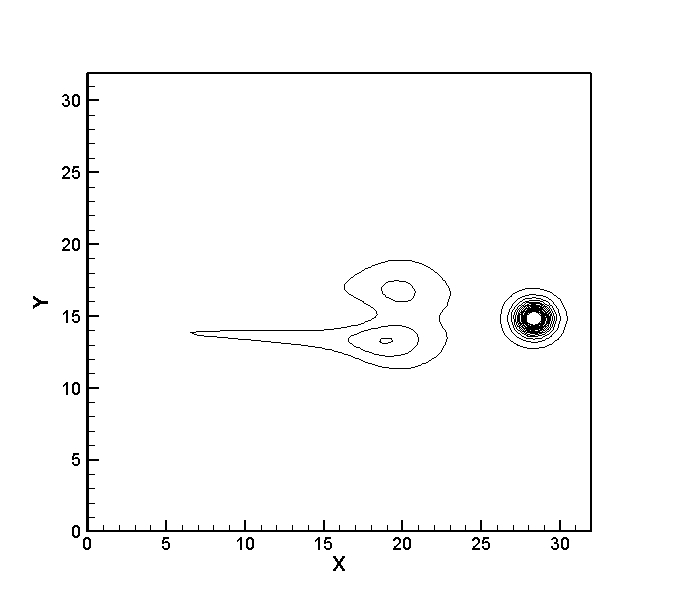}}
    }
    \caption{(continued)}
    \label{fig:edge11}
    \end{center}
    \end{figure}

\section{Conclusions}
\label{sec5}
\setcounter{equation}{0}
\setcounter{figure}{0}
\setcounter{table}{0}

In the previous sections we have studied
a hybrid LDG-HWENO scheme for solving KdV-type equations.
The scheme uses the cell averages of the physical solution and first moment(s) as unknown variables
(a feature of HWENO) while approximates high order spatial derivatives using the local DG method.
It has less unknown variables than a pure LDG method and can be applied to problems involving
high order spatial derivatives. Numerical results have been presented for a selection of one and two dimensional
linear and nonlinear examples. They confirm the designed convergence order  of the scheme.

The new scheme employs an explicit TVD Runge-Kutta method for time integration
and is subject to the CFL condition which requires an extremely small time step
restriction of the form $\Delta t\le C\Delta x^3$.
How to avoid this restriction using suitable implicit time discretization will be the subject of
future work.

\end{document}